\documentclass{amsart}
\usepackage{amsmath,amsfonts,amsthm,amscd,amssymb,pstricks,graphicx,subfig,sidecap}

\newtheorem{thm}{Theorem}[section]
\newtheorem{pro}[thm]{Proposition}
\newtheorem{cor}[thm]{Corollary}

\newtheorem{conjecture}[thm]{Conjecture}

\newtheorem{defn}[thm]{Definition}
\newcommand{\dime}{\operatorname{dim}}

\begin{document}

\author{Jonathan Pakianathan and Troy Winfree}
\title{Quota Complexes, Persistent Homology and the Goldbach Conjecture}

\begin{abstract}
In this paper we introduce the concept of a quota complex and study how the topology of 
these quota complexes changes as the quota is changed. This problem is a simple ``linear" version 
of the general question in Morse Theory of how the topology of a space varies with a parameter.  
We give examples of natural and basic quota complexes 
where this problem frames questions about the distribution of primes, squares and divisors in number 
theory and as an example provide natural topological formulations of the prime number theorem, 
the twin prime conjecture, Goldbach's conjecture, Lehmer's Conjecture, the Riemann Hypothesis and the 
existence of odd perfect numbers among other things.

We also consider random quota complexes associated to sequences of independent random variables 
and show that various formulas for expected topological quantities give L-series and Euler product 
analogs of interest.


\end{abstract}

\maketitle

\tableofcontents

\section{Introduction}
In the area of voting theory in political science, monotone yes/no voting systems are studied. 
If $V$ is a set of voters, then a subset of voters, 
$C \subset V$ is called a losing coalition if the voters in $C$ 
are not sufficient to force an initiative to pass. One can form a simplicial complex 
(see \cite{H} or \cite{Mu} for background on simplicial complexes)
$X$ with vertex set $V$ to encode the voting system by declaring $F=[v_0, \dots , v_n]$ to be a face 
of $X$ if and only if $\{v_0, \dots, v_n\}$ is a losing coalition. The monotone property of the voting system 
then guarantees that any face of a face in $X$ is also a face in $X$ and so $X$ is indeed a 
simplicial complex. In voting theory, it is shown that every voting system can be weighted so 
that it is a quota system. In this paper we study the topological behaviour of various quota 
complexes as the quota is changed. The issues involved are essentially those involved in Morse theory 
or persistant homology where the change in topology of a set is studied as a parameter is increased.

For example in \cite{EH} and \cite{Ka}, the topology of 
a finite union of balls of radius $r$ in $\mathbb{R}^n$, is studied as a function of $r$ for 
fixed centers. This is used both in generating random complexes and in studying the persistant 
shapes/homology of random data sets.

In this paper we first discuss some basic topological properties of quota sets. We then 
illustrate the theory with examples in arithmetic involving the distribution of primes, squares, cubes 
and divisors of a fixed number.

In all these cases, the quota complexes encode significant distributional information and 
our aim is to illustrate some of this.

Let $\mathbb{R}_+$ be the set of positive real numbers.

\begin{defn} 
Let $V$ be a vertex set. A scalar-valued quota system on $V$ is 
given by a weight function $w: V \to \mathbb{R}_+$ and quota $q > 0$.
The quota complex $X[w: q]$ is the simplicial complex on vertex set $V$ 
such that a face $F=[v_0,\dots,v_n]$ is in $X[w:q]$ if and only if 
$w(F)=\sum_{i=0}^n w(v_i) < q$.

\end{defn}

For example, if we have 3 vertices $\{a,b,c\}$ with weights $\{2, 3, 5\}$ respectively then 
at quota $q=9$, the corresponding quota complex would be the boundary of a triangle 
with corners $\{a, b, c\}$. The interior would not be included as the face $[a,b,c]$ has 
weight $2+3+5 = 10 > 9$. If the quota is raised to $11$, this face would be included 
to give a solid triangle whereas if the quota were lowered to $8$, the edge $[b,c]$ would 
be excluded from the complex as it has weight $3+5=8$.

\begin{defn}
Let $V$ be a vertex set. A vector-valued quota system on $V$ is 
given by a weight function $\hat{w}: V \to \mathbb{R}_+^s$ and 
quota $\hat{q} \in \mathbb{R}_+^s$.
The quota complex $X[\hat{w}:\hat{q}]$ is the simplicial complex on vertex set 
$V$ such that a face $F=[v_0,\dots,v_n]$ is in $X[\hat{w}:\hat{q}]$ if and only if 
$\sum_{i=0}^n w_j(v_i) < q_j $ for {\bf some} $1 \leq j \leq s$.

It is easy to see $X[\hat{w}:\hat{q}] = \bigcup_{j=1}^s X[w_j:q_j]$.
$s$ is refered to as the weight dimension of the quota system.

\end{defn}

In the appendix, it is shown that every finite simplicial complex is isomorphic to a vector-valued quota 
complex $X[\hat{w} : \hat{q}]$ where the weights and quota can be taken to be vectors with positive 
integer entries. 

In section~\ref{section: topology}, the basic topology of quota complexes is studied 
especially scalar-valued quota complexes. One of the key results is (all the relevant 
topological definitions can be found in that section):

\begin{thm}
Let $X=X[w: q]$ be a scalar valued quota complex, then $X$ is homotopy equivalent to a 
bouquet of spheres. Let $v_0$ be a vertex of minimal weight, then there is one 
sphere of dimension $s$ in the bouquet for every face $F$ of dimension $s$ in $X$, not containing $v_0$,
such that $q-w(v_0) \leq w(F) < q$.

If $X$ is a vector valued quota complex of weight dimension $N$ then $X$ can be covered 
by $N$ scalar valued quota complexes. As long as $X$ has no "shell vertices" we have 
$Cat(X) \leq 2N-1$ where $Cat(X)$ denotes 
the category of the space $X$. Thus $\frac{Cat(X)+1}{2}$ provides a homotopy invariant 
lower bound on the weight dimension of any quota complex for $X$ with no shell 
vertices.
\end{thm}

These quota systems are then applied to examples in number theory (for the basic 
background needed in number theory see for example \cite{A} or \cite{Te}) :

The prime complex is the full simplicial complex with vertex set equal to the set of primes 
$P=\{2,3,5,7,11,\dots\}$. For a fixed integer $q \geq 3$, $Prime(q)$ is the quota complex 
on vertex set $P$ and quota $q$ (thus the actual vertex set will be all primes less than $q$ 
and the faces will be collections of such primes whose sum is less than $q$.)

In section~\ref{section: prime complex} we study this basic prime complex and show

\begin{thm}
Let $q=2k$ be an even integer $\geq 4$, then $[q-2,q) \cap \mathbb{Z} = \{ O, E \}$ 
where $O=2k-1$, $E=2k-2$ are the unique odd (respectively even) integer in the 
interval $[q-2,q)$. Then the prime complex $Prime(q)$ has the homotopy type of a 
bouquet of spheres where there is one $j$-sphere in the bouquet for every way of writing 
an element in $\{O, E\}$ as a sum of $j+1$ distinct odd primes.

The dimension of the reduced $j$th homology of $Prime(q)$, 
$\dime{\bar{H}_j(Prime(q),\mathbb{Q} )}$ is equal to the number of ways of writing 
an element in $\{O, E\}$ as a sum of $j+1$ distinct odd primes.

For $q \geq 6$, $Prime(q)$ is not connected if and only if $O=q-1$ is a prime number.

For $q \geq 6$, $Prime(q)$ has a non simply-connected component 
if and only if $E=q-2$ is a sum of two distinct odd primes.

Thus the Twin Prime Conjecture is equivalent to the statement that $Prime(q)$ and $Prime(q+2)$ 
are both disconnected complexes for infinitely many values of $q$.

The Goldbach Conjecture is equivalent to the statement that $Prime(q)$ has a non simply connected 
component 
for all $q \geq 6$, with $q$ not equal to twice an odd prime. (The Goldbach conjecture trivially 
holds for $q=2p$, $p$ an odd prime anyway)

\end{thm}

Thus while the $Prime(q)$ complex is defined as a natural quota complex with vertices 
the set of primes less than $q$ and faces given by collections of such primes whose sum is 
less than $q$, the topology of the complex $Prime(q)$ is carried by the sums whose sum 
lies in the shell $[q,q-2)$. Hence one can view the statement above as a topological form of the 
sieve method important in number theory.

In the section~\ref{section: prime complex} we also include data that shows how 
$H_i(Prime(q))$ varies as a function of $q$ for fixed $i$. The behaviour observed is similar to 
that found in the study of random simplicial complexes (see \cite{Ka}) and to behaviour 
observed in the birth-death process in the theory of continuous Markov chains (see \cite{LR}).

In section~\ref{section: Euler} we study how another topological quantity namely the 
Euler characteristic of $Prime(q)$ varies as the quota $q$ is changed. We find 
$$
\chi(Prime(q)) = -\sum_{n=2}^{\infty} \mu(n) L_q(n)
$$
where $\mu$ is the M\"obius function and $L_q$ is a certain characteristic function which converges pointwise to the characteristic 
function of the square-free integers as $q \to \infty$. In order to clarify the issue, we introduce 
the LogPrime quota complex which is the simplicial complex on the set of primes as vertices but with 
the weight of a vertex $p$ being $ln(p)$. We then study how the Euler characteristic of this 
quota complex changes with quota and obtain:

\begin{thm} Let $LogPrime(q)$ be the LogPrime complex with quota $q > 2$. Then
$$
\chi(LogPrime(q)) = - \sum_{2 \leq n < e^q} \mu(n)
$$
where $\mu$ is the M\"obius function. The Riemann hypothesis, that the nontrivial zeros 
of the Riemann zeta function lie on the critical line, is equivalent to the growth condition 
$|\chi(LogPrime(q))| = O(e^{q(0.5 + \epsilon)})$ for all $\epsilon > 0$. In fact, in this case 
the reciprocal of the zeta function will equal a certain difference of $L$-functions associated 
to $\chi(LogPrime(q))$ for $Re(s) > \frac{1}{2}$.
\end{thm}

The substance of this theorem is really due to the work of Titchmarsh (see \cite{Ti}) who gave 
an equivalent statement to the Riemann Hypothesis based on the rate of growth of Merten's function.
However we include it to point out that the problem also can be cast as a question about 
quota complexes.

In section~\ref{section: square complex} we give data for similar complexes encoding 
the distribution of integer squares and cubes. For example one can let 
$Square(q)$ be the simplicial complex on vertices the positive integer squares less than $q$ and 
with the faces consisting of collections of positive integer squares whose sum is less than $q$.
Similarly one can define a complex $Cube(q)$ where cubes replace squares.

One then gets similar flavor theorems:

\begin{thm}
If $Square(q)$ is the square complex, then $Square(q)$ is homotopy equivalent to 
a bouquet of spheres and 
$\dime(\bar{H}_j(Square(q)))$ is equal to 
the number of ways to write $q-1$ as a sum of $j+1$ distinct positive integer squares $> 1$.
Thus $Square(q)$, $q \geq 3$, is connected if and only if $q-1$ is not a positive integer square. 
It is simply connected if and only if $q-1$ is not a positive integer square or the sum 
of two distinct positive integer squares $ > 1$.

If $Cube(q)$ is the cube complex, then $\dime(\bar{H}_j(Cube(q)))$ is equal to the 
number of ways to write $q-1$ as a sum of $j+1$ distinct positive integer cubes $>1$.

Thus the change in the homology of these complexes as quota is varied encodes various 
Waring type problems.

\end{thm}

Data describing the growth of $H_i$ for these two complexes as a function of $q$ is also 
presented in that section.

As a final arithmetical application, in section~\ref{section: divisor complex} we define for 
any integer $n \geq 2$, the divisor complex $Div(n)$ whose vertices consist of the proper 
positive integer divisors of $n$ and whose faces consist of collections of such divisors 
whose sum is less than $n$. (Hence we are using quota $n$).
This complex encodes the distribution of the divisors of $n$, keeping track of which collections 
of divisors have sums less than $n$ or greater or equal to $n$.

Recall an integer $n \geq 2$ is called deficient if the sum of its proper divisors is less than $n$, 
perfect if the sum of its proper divisors is equal to $n$, and abundant if the the sum of its proper divisors is greater than $n$. We obtain the following theorem:

\begin{thm} For any $n \geq 2$, $Div(n)$ is homotopy equivalent to a bouquet of spheres 
where there is one $j$-sphere for every collection of $j+1$ proper divisors of $n$, not including 1 
which sum to $n-1$.

Thus if $n$ is deficient, $Div(n)$ is contractible.

Thus $n$ is perfect if and only if $Div(n)$ is homotopic to a sphere of dimension $\tau(n)-3$ where 
$\tau(n)$ is the number of positive integer divisors of $n$.

\end{thm}

The complex $Div(n)$ in the case that $n$ is abundant can be relatively complicated. Data 
is presented in section~\ref{section: divisor complex} which shows an example of an odd number whose 
complex $Div(n)$ is spherical and of dimension close to $\tau(n)-3$ but not equal. 
This odd number is ``close" to being an odd perfect number in a topological sense; of course 
the existence of an actual odd perfect number is still open. Even perfect numbers are in bijective 
correspondence with Mersenne primes and it is an open question whether 
there are infinitely many of these.

In section~\ref{section: counting}, a generating function for Euler characteristics of quota 
complexes is found and a quota-complex formulation of Lehmer's conjecture is provided among other things.

In the final section, we consider a finite set of independent 
continuous random variables $X_1, \dots X_N$ with continuous density functions 
$f_i$ with compact support in $[m,\infty)$ where $m > 0$. We examine 
the random quota complex associated to these random variables with fixed minimal weight 
$X_0=m$ and quota  $q > m > 0$. We derive various 
formulas for the expected topology of this complex (see the section for the relevant 
definitions and \cite{B} for basic background.)

\begin{thm}
Let $X_0=m > 0$. Let $X_1, \dots, X_N$ be independent, continuous random variables 
with density functions $f_1, \dots, f_N$ which are continuous with 
compact support in $[m, \infty)$ and let $\mathbb{X}[q]$ be the random scalar quota 
complex determined by this collection and quota $q > m > 0$.

Then for $j \geq 1$, 
$$
E[\dim(\bar{H}_{j-1}(\mathbb{X}[q], \mathbb{Q}))] = \sum_{\mathfrak{J}, |\mathfrak{J}|=j} (f_{\mathfrak{J}} \star \mathbb{I}_m)(q)
$$
is a continuous function of $q$ with compact support, where $\star$ denotes convolution.

Furthermore we have 
$$
1 - E[\chi(\mathbb{X}[q])] = \frac{1}{4 \pi^2 i} \int_{-\infty}^{\infty} e^{2 \pi i \alpha x} \left(\prod_{j=0}^N (1-\hat{f}_j(\alpha)) - (1-\hat{f}_0(\alpha))\right) \frac{d\alpha}{\alpha}.
$$
is a continuous function of $q$ with compact support, where $\hat{f}$ denotes the 
Fourier transform of $f$.
\end{thm}

We then give an example where the final equality in the last theorem is:
$$
\sum_{1 \leq n < e^q} \mu(n) = \frac{1}{4 \pi^2 i} \int_{-i\infty}^{+i\infty} e^{sx}\left(\prod_{p \in P}(1-\frac{1}{p^s}) - (1-\frac{1}{2^s})\right) 
\frac{ds}{s}
$$
where $P$ is the set of primes less than $e^q$. Note for $Re(s)>1$, 
$\prod_{p \in P} (1-\frac{1}{p^s}) \to \frac{1}{\zeta(s)}$ as $q \to \infty$ where $\zeta(s)$ 
is the Riemann zeta function.

\section{Topology of quota complexes}
\label{section: topology}

In this section we discuss the topological results for scalar quota complexes. For the basic
background needed in this section see \cite{H} or \cite{Mu}. The reader more interested in applications to 
number theoretic examples can skim the statements 
of theorems and move on to the application sections of the paper.

First we recall some elementary definitions for simplicial complexes.

\begin{defn} Let $X$ be a simplicial complex and let $v$ be a vertex in $X$.
The closed star of $v$, $\bar{St}(v)$ is the union of all faces in $X$ that contain $v$.
The open star of $v$, $St(v)$ is the union of the interior of all faces in $X$ that contain $v$.
The link of $v$ is defined as $Lk(v)=\bar{St}(v)-St(v)$. Thus $Lk(v)$ consists of all faces 
$F=[v_0,\dots,v_k]$ in $X$ not containing $v$ such that $[v_0,\dots,v_k,v]$ is also a face of $X$.

The open star $St(v)$ is an open neighborhood of $v$ in $X$, the closed star $\bar{St}(v)$ is its 
closure and $Lk(v)$ is its boundary in $X$.
\end{defn}

The topological structure of a scalar valued quota complex is controlled to a large degree 
by the part of the complex outside the closed star of a vertex of minimal weight as the next theorem 
illustrates:

\begin{pro}
\label{pro: scalarquotaprop}
Let $X$ be a scalar weighted finite simplicial complex. Let $v_{min}$ be a vertex of minimal 
weight. If $F$ is a face of $X$ that is not in $\bar{St}(v_{min})$ then 
the boundary of $F$ is completely contained in $Lk(v_{min}) \subset \bar{St}(v_{min})$.

Furthermore $F$ is such a face if and only if $F$ does not contain $v_{min}$ and 
$q - w(v_{min}) \leq w(F) < q$.

\end{pro}
\begin{proof}
 If $V$ is the vertex set of $X$, $w: V \to \mathbb{R}_+$ the weight function and $q > 0$ the quota then 
$[v_0,\dots, v_k]$ is a $k$-face of $X$ if and only if $\sum_{i=0}^k w(v_i) < q$. 
Let $v_{min}$ be a vertex of minimal weight which exists as $V$ is finite.

Note that if $F=[v_0,\dots,v_k]$ is a face of $X$ not contained in the closed star of $v_{min}$ 
then $v_{min} \notin F$ and $w(F) < q$ as $F$ is a face of $X$. Furthermore,  
$[v_0,\dots,v_k,v_{min}]$ cannot be a face in $X$ as if not it would be in $\bar{St}(v_{min})$ 
and hence $F$ being a face of $[v_0,\dots,v_k,v_{min}]$ would also be in 
$\bar{St}(v_{min})$. Thus $w(F) + w(v_{min}) \geq q$. Putting the inequalities together gives 
$q-w(v_{min}) \leq w(F) < q$. Conversely, reversing the argument, it is easily checked 
that any face $F$ of $X$ not containing $v_{min}$ and 
satisfying these final inequalities is not in the closed star of $v_{min}$.

Let $F=[v_0, \dots v_k]$ and let $\sigma = [v_0, \dots, \hat{v}_i, \dots, v_k]$ be a face in the boundary 
of $F$ obtained by removing the $i$th vertex. Note that $w(\sigma)=w(F)-w(v_i) \leq w(F)-w(v_{min})$ 
and so $w([\sigma,v_{\min}])=w(\sigma) + w(v_{min}) \leq w(F) < q$ and so 
$[\sigma,v_{min}]$ is a face in $X$ and so $\sigma \subseteq Lk(v_{min}) \subset \bar{St}(v_{min})$
as desired.

Thus the boundary of $F$ is completely contained in $Lk(v_{min}) \subset 
\bar{St}(v_{min})$ as desired.

\end{proof}

For the next theorem, we recall the basic fact that if $A$ is a contractible subcomplex 
of a CW or simplicial complex $X$ then the quotient map $\pi: X \to X/A$ which collapses 
$A$ to a point is a homotopy equivalence. However even if $A$ and $X$ are simplicial 
complexes, the collapsed space $X/A$ is only a $CW$-complex in general.

We will apply this fact to the case of a finite scalar quota complex $X$ where $A$ will be 
the closed star of a minimal weight vertex. Closed and open stars of a vertex $v$ 
are always contractible as they are star-convex with respect to the vertex $v$.

In this context, it is important to note that if $F=[v_0,\dots,v_k]$ is a $k$-face of $X$ 
then the space $F/\partial{F}$ obtained by collapsing the boundary of $F$ to a point is 
homeomorphic to a  $k$-sphere $S^k$ i.e., the space of unit vectors in $\mathbb{R}^{k+1}$.
Recall that a bouquet of spheres is the wedge product of a collection of spheres (not necessarily 
of the same dimension and the $0$-sphere is allowed). Intuitively this is a collection of 
spheres attached at a common point.

\begin{thm}[Scalar quota complexes are homotopy equivalent to bouquets of spheres]
\label{thm: scalarquotathm}
Let $X$ be a finite scalar quota complex and let $A=\bar{St}(v_{min})$ be the 
closed star of a vertex of minimal weight. Then the quotient map 
$\pi: X \to X/A$ is a homotopy equivalence and 
$X/A$ is a bouquet of spheres where there is one $i$-sphere for each 
$i$-face $F=[v_0, \dots, v_i]$ in $X$ not containing $v_{min}$ and such that 
$q-w(v_{min}) \leq w(F) < q$.

Thus the reduced integer homology groups of $X$, $\bar{H}_i(X)$ are free abelian groups of finite rank 
equal to the number of $i$-faces $F$ of $X$, not containing $v_{min}$ 
with $q-w(v_{min}) \leq w(F) < q$.
\end{thm}
\begin{proof}
The fact that the collapse map $\pi: X \to X/A$ is a homotopy equivalence was explained 
in the paragraphs before the statement of the theorem. Any face of $X$ inside 
$A=\bar{St}(v_{min})$ maps to the collapse basepoint $A/A$ in $X/A$. 

By proposition~\ref{pro: scalarquotaprop},  the $i$-faces $F$ of $X$ not contained in $A$ 
are exactly the $i$-faces $F$ not containing $v_{min}$ with 
$q - w(v_{min}) \leq w(F) < q$ and such a face maps to an $i$-sphere $F/\partial F = S^i$ 
in $X/A$. From this it is easy to see that $X/A$ is a bouquet of spheres as claimed. 

The comment on reduced integer homology follows immediately from basic facts 
about homology.
\end{proof}

Note Theorem~\ref{thm: scalarquotathm} is a sort of topological sieve in the sense that 
it says the interesting topology of a scalar quota complex is carried by the "shell faces" 
i.e. the faces $F$ with $q-w(v_{min}) \leq w(F) < q$ so that their weight is concentrated 
in a narrow shell near the quota $q$. 

As we will not use vector-weighted quota complexes in the applications of the rest of the paper, 
the theorems and proofs about the vector-weighted case are deferred to Appendix~\ref{section: vector}.
The reader is also referred to Appendix~\ref{section: Alves}, for the proof that every finite simplicial complex is isomorphic to a (vector weighted) quota complex. Thus for example, every closed manifold is homeomorphic 
to some vector weighted quota complex. The picture for vector weighted quota complexes is much less complete than that for scalar weighted quota complexes. In the applications 
in the rest of the paper, scalar weighted quota complexes are used predominantly because 
their topological structure is completely determined up to homotopy by the results of this section.

\section{Application: The Prime Complex and the Goldbach conjecture}
\label{section: prime complex}

In the current and following sections we examine three examples of scalar quota complexes that
 arise from sequences of increasing positive integers $V = \{v_n\}_{n=1}^{\infty}$. In general, for 
 integers $q> v_1$  we take  $V(q)$ to be the quota complex on the vertex set $V$ with quota $q$. 
 As described in section~\ref{section: topology} the topology of $V(q)$ is entirely determined by the 
 number of ways to add integers $v_n\in V$, where $v_1 \ne v_n < q$, so that the sum falls in the 
 interval $[q-v_1, q)$. So $V(q)$ topologically encodes the number of ways to express integers in 
 $[q-v_1, q)$ as sums of distinct elements in $V-\{v_1\}$. Our main interest is in describing the 
 behavior of $V(q)$ as $q$ is increased for the cases where $V$ is the set of primes, squares and cubes.

Our description of $V(q)$ will be data-based and will center on the functions $h_i(q)$ and $s_i(q)$, 
where $h_i(q) = \dim (\bar{H}_i(V(q)),\, \mathbb{Q})$ and $s_i(q)$ is the number of $i$-simplexes in
$V(q)$ not containing $v_1$. So $s_i(q) = \left| \left\{ \, U \subset V - \{ v_1 \} \, : \,  |U| = i+1,\, \sum_{v\in U} v < q\,  \right\}\right|$ and
 $h_i(q) = s_i(q) - s_i(q-v_1)$ by Theorem~\ref{thm: scalarquotathm}. Furthermore note that for $i>0$, $s_i(q)+s_{i-1}(q-v_1)$ is the number of 
 $i$-simplexes in $V(q)$ and of course $s_0(q)+1$ is the number of zero simplexes. In order to capture the relative growth of these values we will examine the ratios:
\[ S_i(q) = \frac{s_i(q)}{\sum_j s_j(q)} \hspace{0.1in}\text{and}\hspace{0.1in}H_i(q) = \frac{h_i(q)}{\sum_j h_j(q)}. \]
Note of course that $0 \leq S_i(q), H_i(q) \leq 1$ count the fraction of the cells (respectively homology) 
of the quota complex $V(q)$ that are concentrated in dimension $i$.

Before describing the data, which was produced with a purpose-built algorithm in C{}\verb!++!, we first develop 
a simple theoretical context.  Let $\kappa > 0$ be an integer and suppose there is a monotonically increasing, 
differentiable function $f:[\kappa,\, \infty)\to\mathbb{R}$, with $\lim_{x\to\infty}f(x) = \infty$ and such that 
$s_0(q) \sim f(q)$, by which we mean that $\lim_{q\to\infty}s_0(q)/f(q) = 1$. Such a function $f$ 
will be called an interpolating function for the vertex count function $s_0(q)$.

Set $\widehat{s}_i(x) =  {f(x) \choose i+1}$, 
where we are thinking of ${x \choose k}$ as the degree $k$ polynomial $x(x-1)(x-2)\cdots(x-k+1)/k!$. So 
certainly $\widehat{s}_i(x) \sim f(x)^{i+1}/(i+1)!$. We make the approximation  
\[\widehat{s}_i(q/(i+1)) \sim {s_0(\lceil q/(i+1) \rceil) \choose i+1}  \le s_i(q) \le {s_0(q) \choose i+1} \sim  \widehat{s}_i(q).\]
This comes from noting that $s_0(q) \choose i+1$ is the number of potential 
$i$-simplices in the complex that one can make out of the $s_0(q)$ vertices other than $v_1$ 
and hence certainly bounds $s_i(q)$ from above. On the other hand, 
${s_0(\lceil q/(i+1) \rceil) \choose i+1}$ is the number of $i$-simplices possible 
that one can make out of vertices other than $v_1$ but with weight below $\frac{q}{i+1}$; such
simplices definitely are below quota and hence definitely count towards $s_i(q)$.

In this setup, various general expectations can be derived for the functions $S_i(q)$ and $H_i(q)$.
The expectations and their proofs are collected in proposition~\ref{prime complex: proposition} 
in Appendix~\ref{section: ExpJust}. It is shown there that $\{ S_i: \mathbb{Z} \to [0,1] \,|\, i=0,1,2, \dots \}$ 
should be a family of unimodal functions of $q$, each with a single local maximum which moves 
to the right as $i$ increases and such that the heights of the maxima decrease towards $0$ as $i \to \infty$. 
All of these expectations are realized in the data. 

In the primes quota complex ($v_1 = 2$) 
the data indicates that the family of functions $H_i(q)$ has the same global behavior as $S_i(q)$, whereas in the 
squares and cubes cases ($v_1 = 1$), $H_i(q)$ has no discernible shape outside of appearing to tend to zero as $q\to\infty$.

We now consider the case $V=P=\{2,\, 3,\, 5,\cdots  \}$, the set of primes,  and its 
quota complex $Prime(q)$. For this vertex set we will denote the corresponding simplex and homology functions 
discussed above with a superscript $P$, so for example $s_i^P(q)$ is the number of $i$-simplexes in $Prime(q)$ not containing 2. 
Note that in this case $s_0^P(q) = \pi(q-1) - 1$, where $\pi$ is the prime number counting function.

%%%%%%Primes: Simplex Graphs

\begin{figure}[ht]
  \subfloat[$S_i^P(q)$, $0<q\le550$, $0\le i\le 6$]
  {\label{primesFig:simplexRaw}\includegraphics[width=0.5\textwidth]{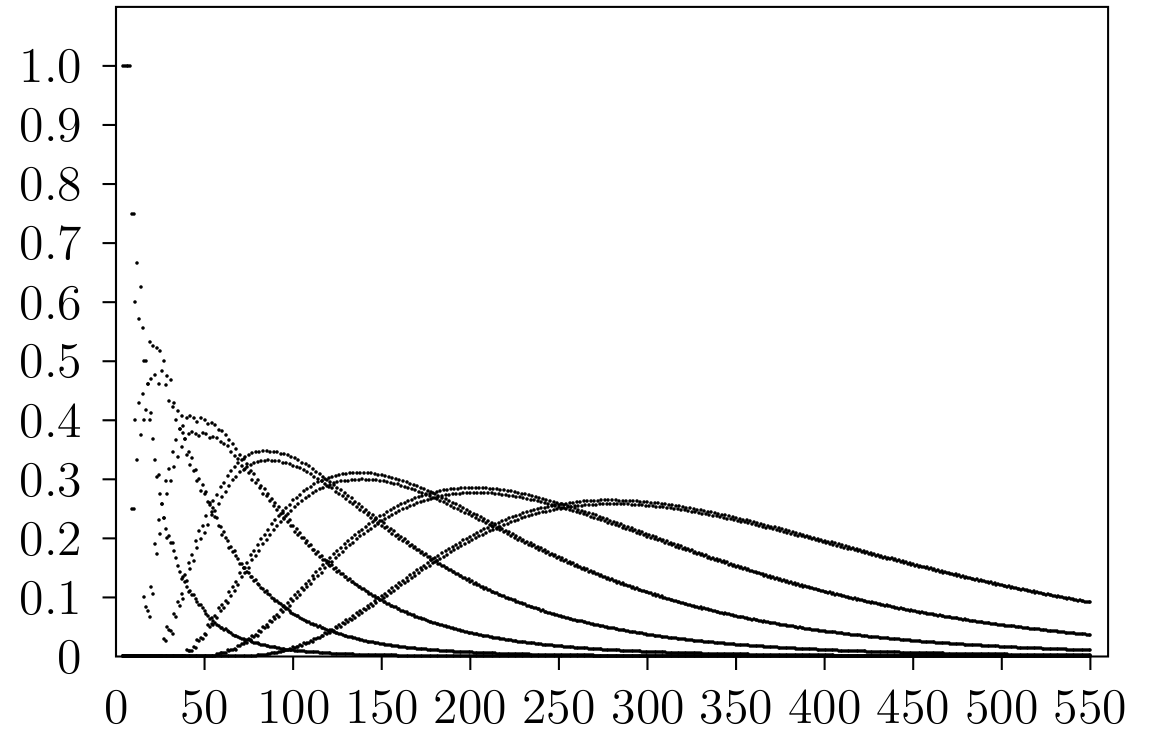}}           
  \subfloat[$S_{i,ave}^{P}(q)$, $0<q\le550$, $0\le i\le 6$]
  {\label{primesFig:smoothSimplex}\includegraphics[width=0.5\textwidth]{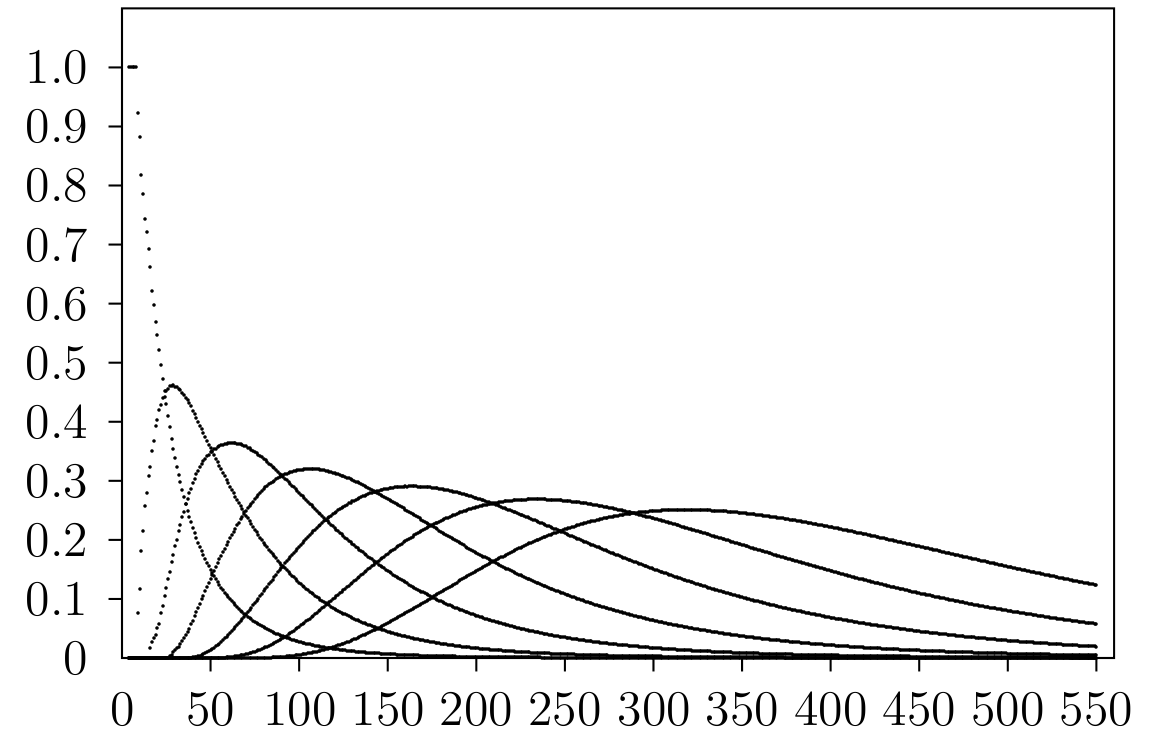}}
  
  \subfloat[$H_i^P(q)$, $0<q\le550$, $0\le i\le 6$]
  {\label{primesFig:homologyRaw}\includegraphics[width=0.5\textwidth]{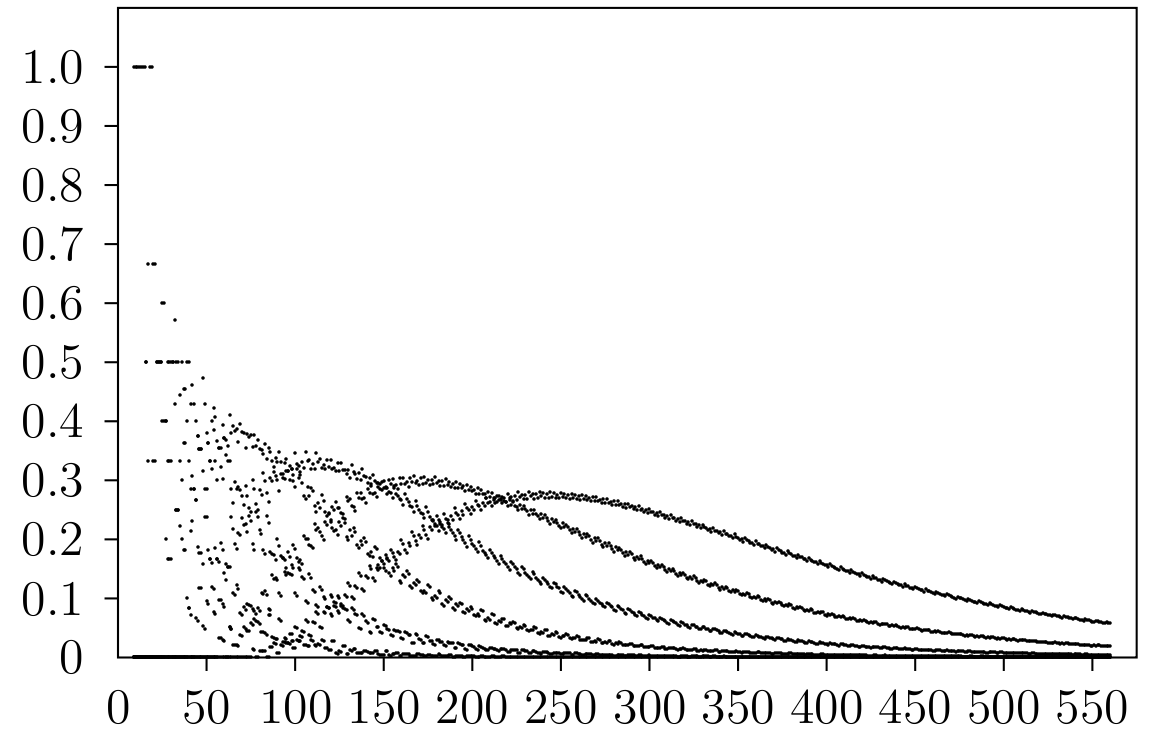}}                
  \subfloat[$H_{i,ave}^{P}(q)$, $0<q\le550$, $0\le i\le 6$]
  {\label{primesFig:smoothHomology}\includegraphics[width=0.5\textwidth]{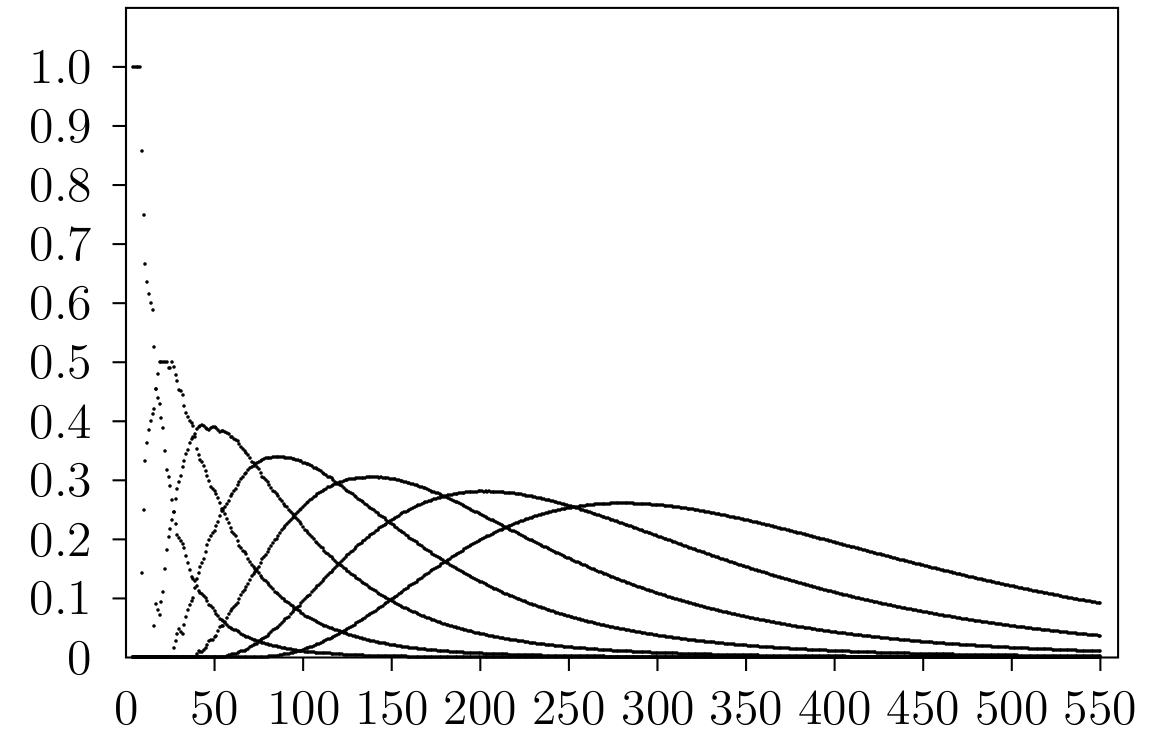}}
  
   \caption{$Prime(q)$}
  \label{primesFig:simplex}
\end{figure}

Recall that the prime number theorem says $\pi(q) \sim f_P(q)$ where $f_P(x) = x/\ln x$. Certainly 
$f_P: [3,\, \infty)\to \mathbb{R}$ has the required properties to be an interpolating function, so we expect $S_i^P(q)$ to have the 
behavior described above, at least for $i\ge3$. Our data was generated using the first 100 odd primes 
(so the primes 3 through 547), producing complete information in the range $0<q\le550$ 
and up to homology dimension $16$.  Figure~\ref{primesFig:simplexRaw} shows $S_i^P(q)$ versus $q$ for 
$0<q\le550$ and cell-dimensions $0\le i \le 6$; here our global expectations are closely exhibited:  $S_i^P(q)$ is 
unimodal and as $i$ increases the maximum height decreases and the position of the maximum moves to the right. 
Figure~\ref{primesFig:smoothSimplex} is a smoothing of $S_i^P(q)$, showing the quantity 
$S_{i,ave}^P(q) = \frac{1}{q}\sum_{k=1}^{q}S_i^P(k)$. The ``double-line'' effect evident in 
figure~\ref{primesFig:simplexRaw} occurs because we are summing odd numbers, so depending 
on the length of the sum only odd or even sums occur in each dimension. 

As described in the introduction, the homology data of $Prime(q)$ encodes several classical 
arithmetic conjectures in number theory. In particular if $h_1^P(q) > 0$ for all $q>8$ then the 
Goldbach conjecture holds, in fact in this case a stronger refinement of Goldbach holds, namely 
that every even number greater than or equal to 8 is the sum of two distinct odd primes. Our data 
is consistent with the stronger refinement of Goldbach and thus we conjecture:

\begin{conjecture}
$\dime(H_1(Prime(q),\mathbb{Q})) \neq 0$ for all $q > 8$.
 
Equivalently, every even number $q \geq 8$ is a sum of exactly two distinct odd primes.
\end{conjecture}

The classical Goldbach conjecture is equivalent to the above conjecture holding for 
all even numbers greater than 3 and 
not equal to twice a prime. (The Goldbach conjecture is trivial for even 
numbers equal to twice a prime anyway.)

Figure~\ref{primesFig:homologyRaw} 
shows $H_i^P(q)$ for $0<q\le550$ and $0\le i\le 6$. Although appearing fairly scattered below 
$q\approx100$ in these dimensions, the same global behavior as $S_i^P(q)$ manifests fairly 
rapidly as $q$ increases. As before, we also provide a smoothing of the data in figure~\ref{primesFig:smoothHomology}.

In general, the quantities $h_i^P(q)$ count the ways to sum to even or odd positive integers using 
$i+1$ distinct odd primes. Specifically if $h_i^P(q) > 0$ for all sufficiently large $q$ and $i+1$ is even/odd, then all sufficiently 
large positive even/odd integers are expressible as a sum of distinct odd primes. At $i=2$ this is a refinement of the weak
Goldbach conjecture that all odd numbers greater than 7 are a sum of three odd primes. The 
weak Goldbach conjecture has been shown to hold for almost all odd numbers by Vinogradav.

In our data $h_2^P(q)>0$ for $20\le q\le550$ and so it is reasonable to conjecture:

\begin{conjecture}
$\dime(H_2(Prime(q),\mathbb{Q})) \neq 0$ for all $q \geq 20$.
 
Equivalently, every odd number $q \geq 19$ is a sum of exactly three distinct odd primes.
\end{conjecture}

%%%%%%Primes: Linear Graph

\begin{figure}[t]
\includegraphics[width=0.5\textwidth]{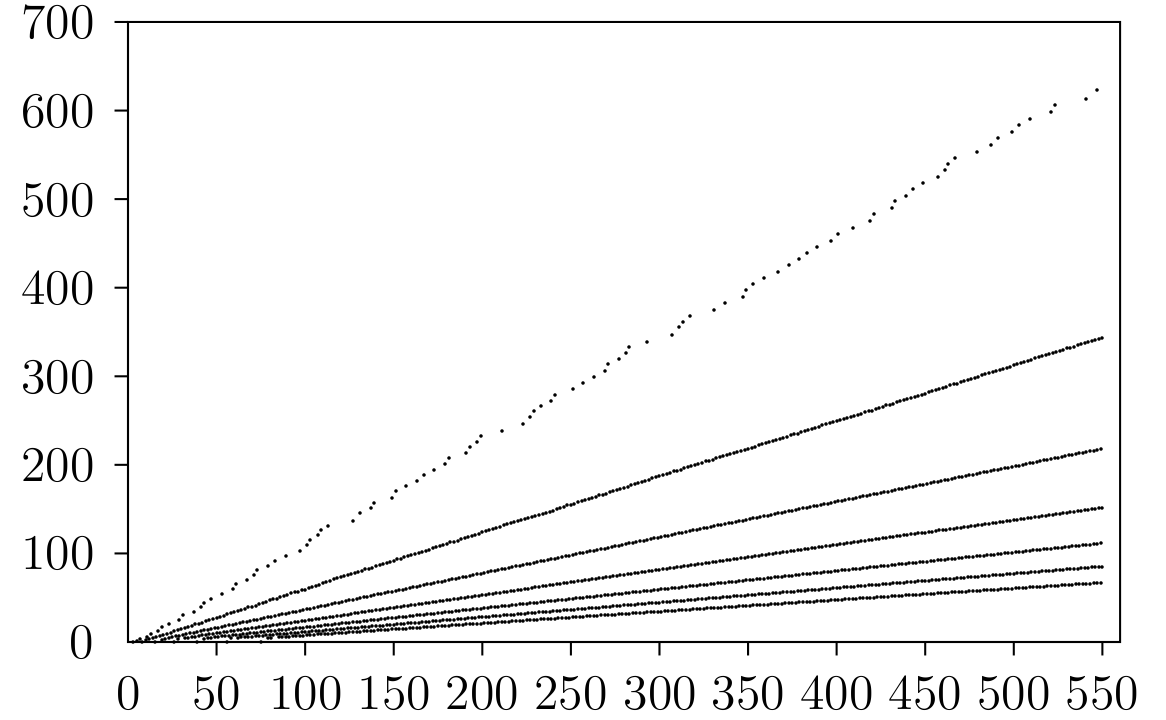}
\caption{$\left(s_i^P(q)\right)^{1/(i+1)}ln(q)$, $0<q\le550$, $0\le i\le 6$}
\label{primesFig:linear}
\end{figure}

Perhaps the most striking feature of the data is displayed in  figure~\ref{primesFig:linear}. Here, 
in a naive attempt to capture the expected asymptotic growth of $s_i^P(q)$, 
we graph $\left(s_i^P(q)\right)^{1/(i+1)} ln(q)$ versus $q$ 
for $0<q\le550$ and $0\le i\le 6$. In this range the
 data displays strong linearity. The data line with slope nearly one is the line for the zero-simplexes, the 
 data line with the next greatest slope is the one-simplexes, and so-on. Note that since 
 $\widehat{s}_i(q) \sim f_P(q)^{i+1}/(i+1)! = \frac{q^{i+1}}{ln(q)^{i+1}(i+1)!}$ by the prime number theorem and $\widehat{s}_i(\frac{q}{i+1}) \leq s_i(q) \leq 
 \widehat{s}_i(q)$, heuristically one might expect lines when plotting $\left(s_i^P(q)\right)^{1/(i+1)}\ln q$ versus $q$ with slopes roughly somewhere between $\frac{((i+1)!)^{-1/(i+1)}}{i+1}$ and $((i+1)!)^{-1/(i+1)}$. Running a least squares approximation on the data gives slopes of $0.632374$ when $i=1$; $0.404613$ when $i=2$; $0.284124$ when $i=3$; $0.211868 $ when $i=4$; $0.164796 $ when $i=5$; and $0.132366$ when $i=6$. These values are consistent with the heuristics. Thus it is 
 reasonable to conjecture:
 
 \begin{conjecture}
 Let $s_i^P(q)$ denote the number of sets of $i+1$ distinct odd primes whose sum is below $q$.
 Then 
 $$
 s_i^P(q) \sim \frac{C_i^{i+1} q^{i+1}}{ln(q)^{i+1}}
 $$
 where 
 $1=C_0 > C_1 > C_2 > \dots > 0$.
 
 Note the case $i=0$ is the prime number theorem which is of course known to be true.
 Approximate values of the constants $C_1,\dots,C_6$ are listed in the paragraph before 
 the conjecture.
 \end{conjecture}
 
\section{Application: Sum of Squares and Cubes complexes}
\label{section: square complex} 

In this section we present data for the quota complexes $V(q)$ with $V = S = \{ 1,\, 4,\, 9,\, 16,\ldots \}$ 
the set of squares, and $V = C = \{1, \, 8,\, 27, \,64,\ldots \}$ the set of cubes, which we denote 
$Square(q)$ and $Cube(q)$ respectively. We indicate the various functions of interest for these 
quota complexes with a superscript $S$ or $C$ as we did in section~\ref{section: prime complex}. 
Note that since $s_0^S(q) = \lfloor\sqrt{q-1}\rfloor - 1$ and $s_0^C(q) = \lfloor\sqrt[3]{q-1}\rfloor - 1$, 
our expectations for the global behavior of the families $S_i^S$ and $S_i^C$ are once again 
informed by proposition \ref{prime complex: proposition}.

%%%%%%Squares: Everything

\begin{figure}[ht]
  \subfloat[$S_i^S(q)$, $0<q\le629$, $0\le i\le 6$]
  {\label{squaresFig:simplexRaw}\includegraphics[width=0.5\textwidth]{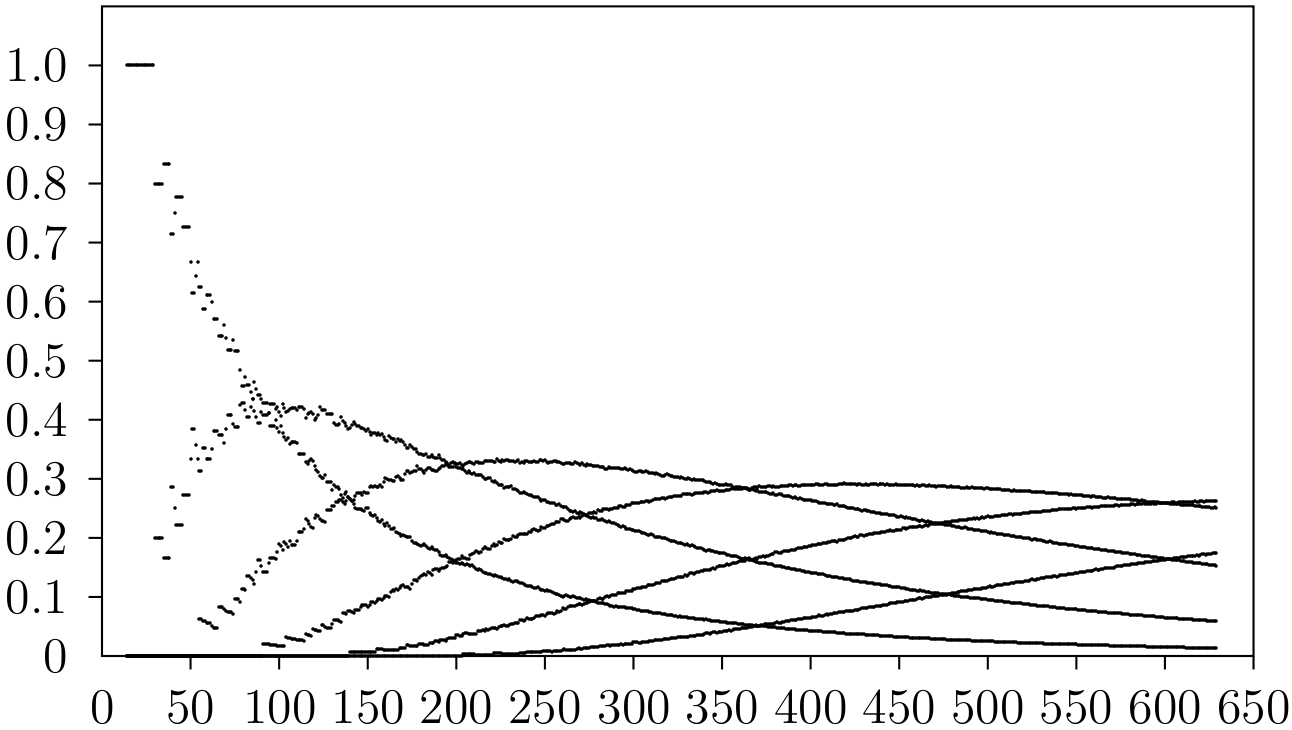}}           
  \subfloat[$S_{i,ave}^S(q)$, $0<q\le629$, $0\le i\le 6$]
  {\label{squaresFig:smoothSimplex}\includegraphics[width=0.5\textwidth]{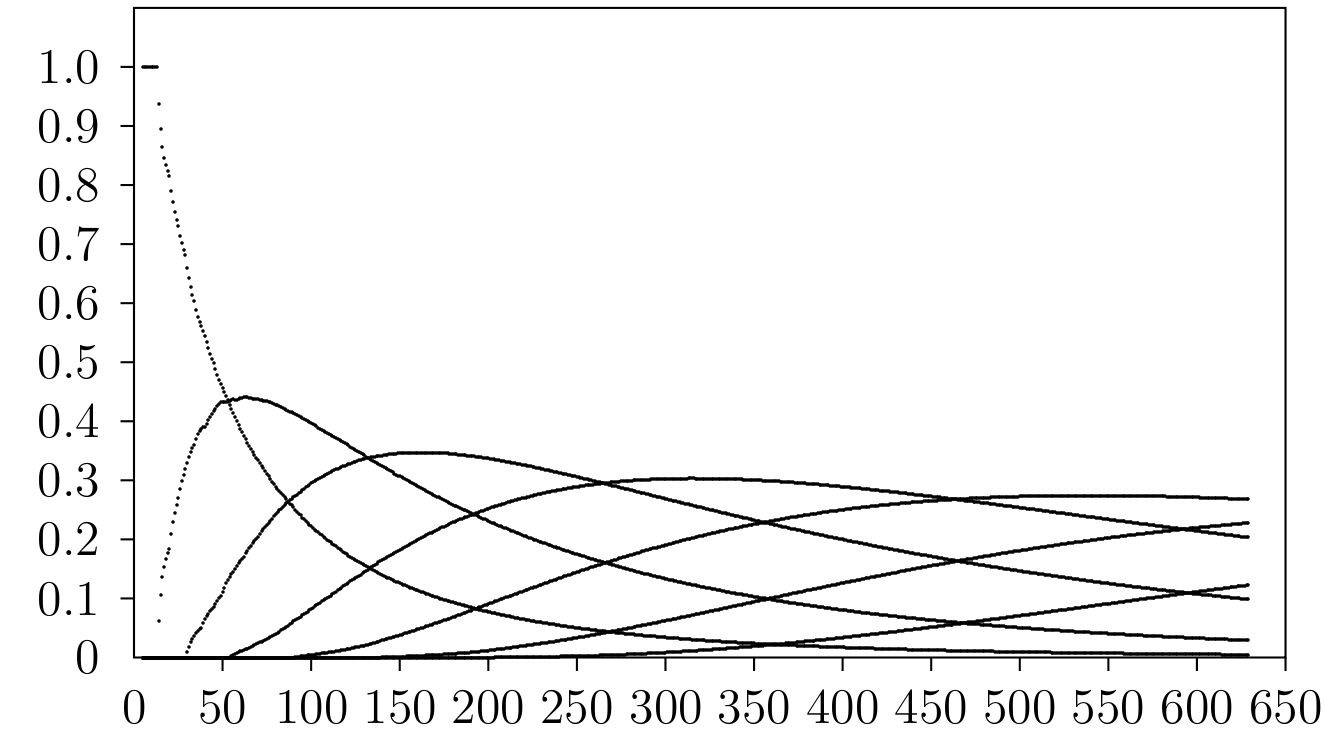}}
  
    \subfloat[$H_i^S(q)$, $0<q\le629$, $0\le i\le 6$]
  {\label{squaresFig:homologyRaw}\includegraphics[width=0.5\textwidth]{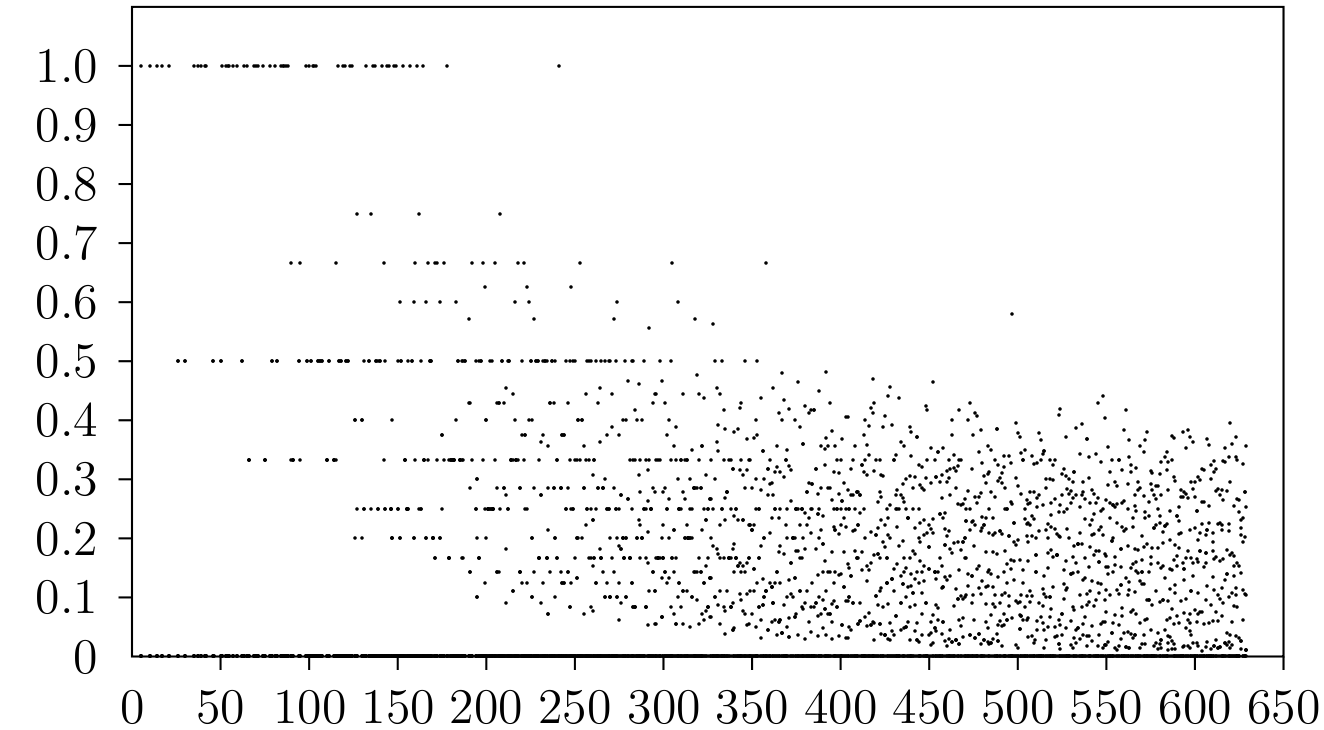}}           
  \subfloat[$s_i^S(q)^{\frac{2}{i+1}}$, $0<q<629$, $0\le i\le 6$]
  {\label{squaresFig:linear}\includegraphics[width=0.5\textwidth]{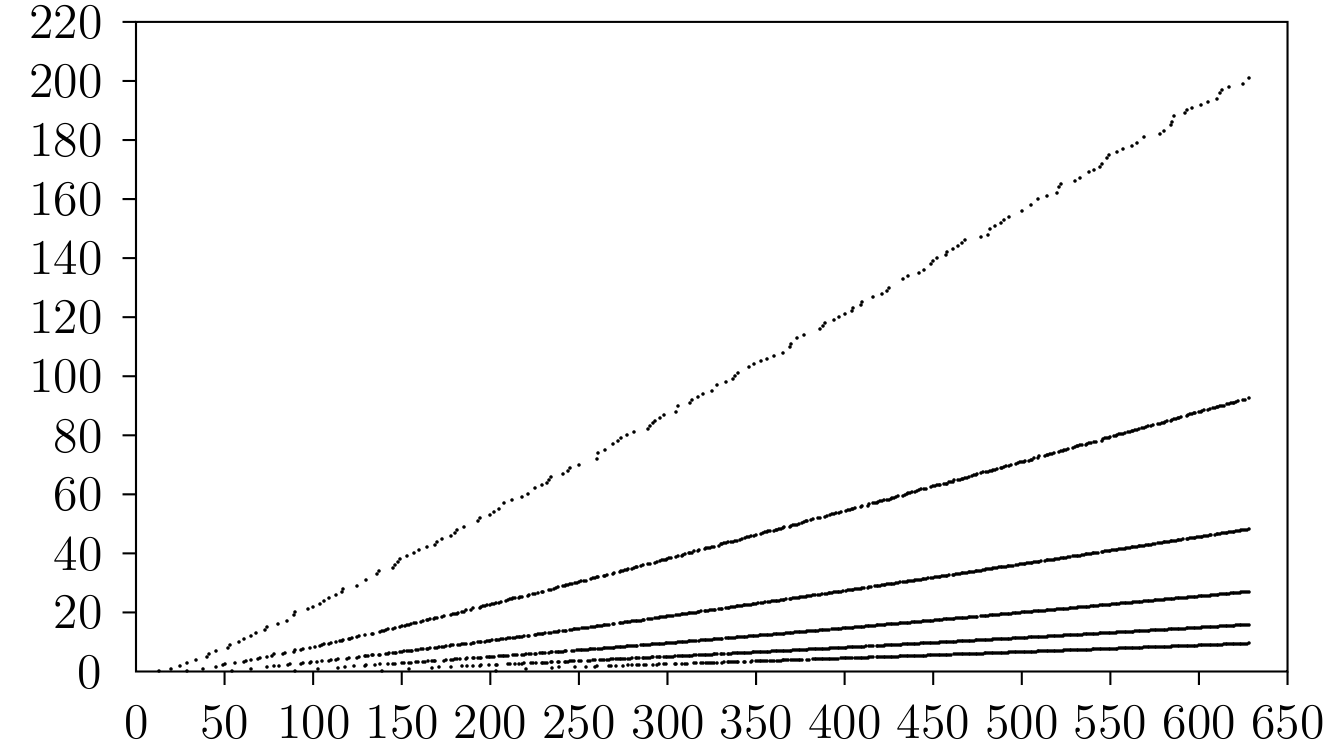}}
  
  \caption{$Squares(q)$}
  \label{figure:squares}
\end{figure}

%%%%%%Cubes: Simplex and linear

\begin{figure}[ht]
  \subfloat[$S_i^C(q)$, $0<q\le15633$, $0\le i\le6 $]
  {\label{cubesFig:simplexRaw}\includegraphics[width=0.5\textwidth]{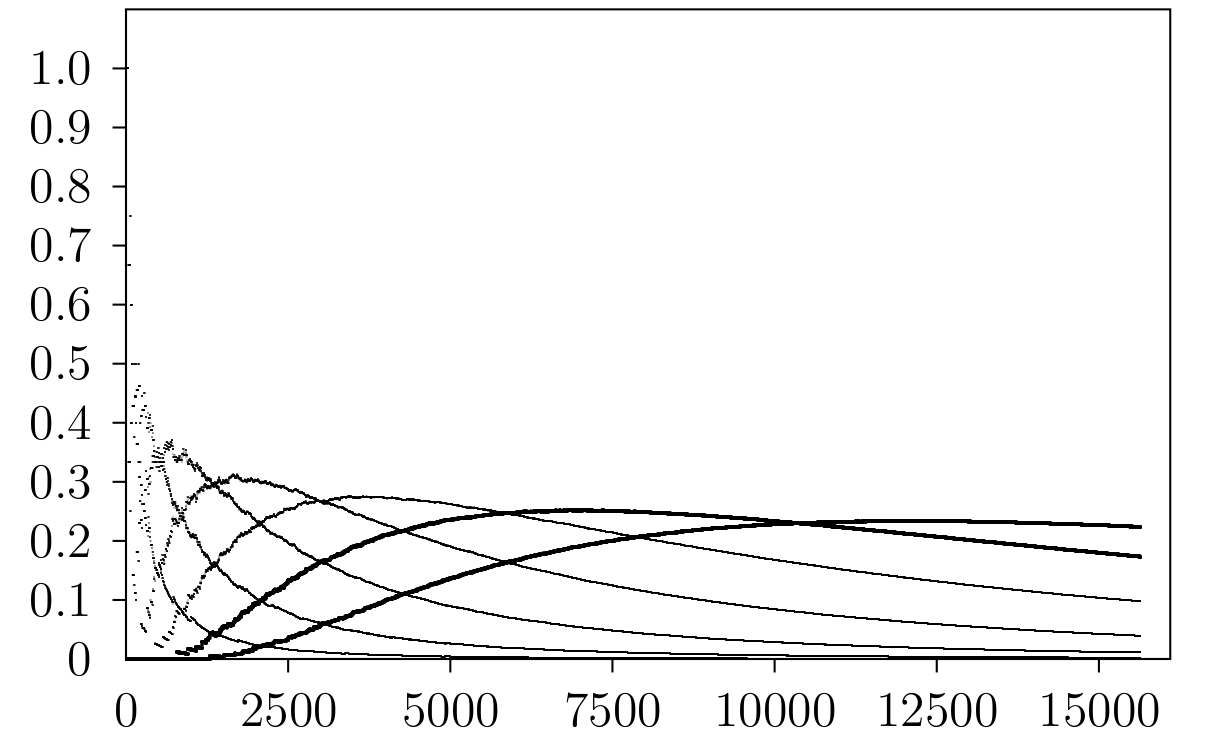}}           
  \subfloat[$s_i^C(q)^{\frac{3}{i+1}}$, $0<q<15633$, $0\le i\le4 $]
  {\label{cubesFig:linear}\includegraphics[width=0.5\textwidth]{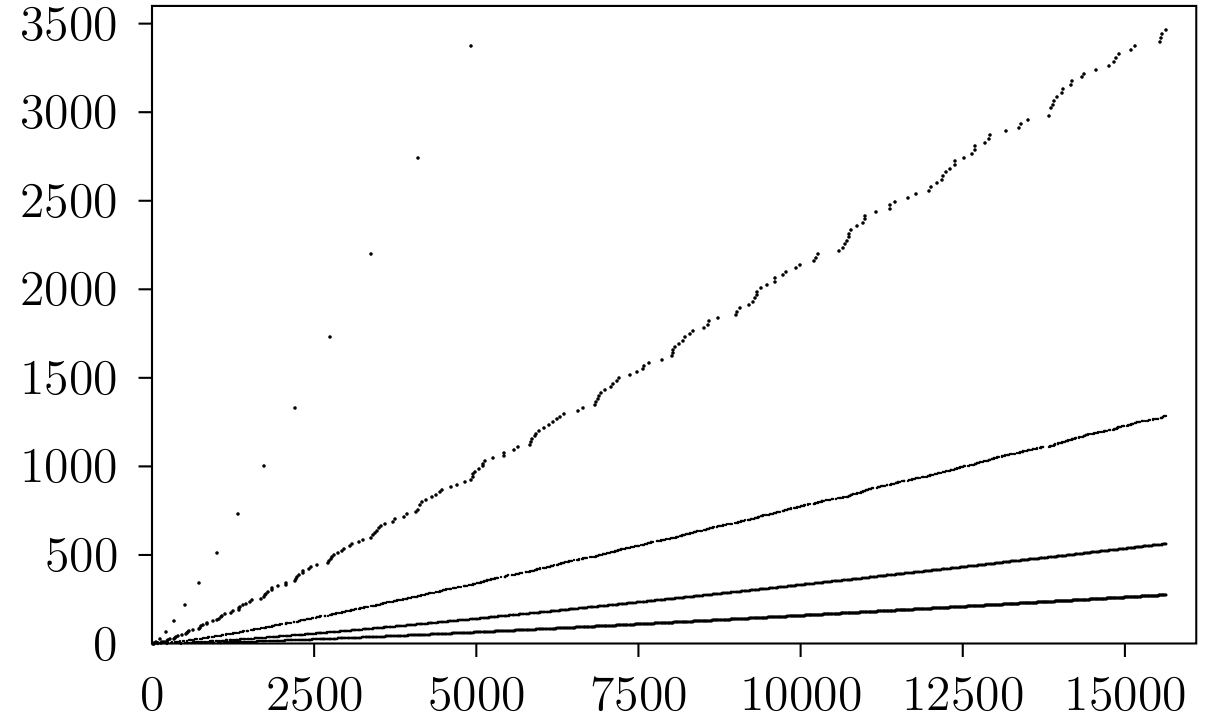}}
  \caption{$Cubes(q)$}
   \label{figure:cubes}
\end{figure}

Data for $Square(q)$ was generated using the first 25 squares, producing complete 
information for $0<q\le 629$ and up to homology dimension 9. Figure \ref{squaresFig:simplexRaw} 
displays $S_i^S(q)$ for $0<q\le 629$ and cell-dimensions $0\le i \le 6$. Here our expectations are strongly evident 
and are made all the more so in figure \ref{squaresFig:smoothSimplex} which shows 
$S_{i,ave}^S(q) = \frac{1}{q}\sum_{k=1}^{q}S_i^S(k)$. The homology data $H_i^S(q)$ is notable 
for its lack of shape and is displayed in figure \ref{squaresFig:homologyRaw}. Finally, the heuristic growth of $s_i^S(q)$ as $q^{\frac{i+1}{2}}$ is verified in figure \ref{squaresFig:linear} as a graph of $s_i^S(q)^{\frac{2}{i+1}}$ versus $q$.

Data for $Cube(q)$ was generated using the first 25 cubes, producing complete information 
in the range $0<q\le15633$ and up to homology dimension 13. Once again our expectations are 
closely realized in the data: figure \ref{cubesFig:simplexRaw} shows $S_i^C(q)$ for $0\le i\le 6$ and 
figure \ref{cubesFig:linear} displays the growth of $s_i^C(q)$ as  $q^{\frac{i+1}{3}}$. A graph of $H_i^C$ is not 
included as this data is scattered and apparently formless.

Based on the data, it is reasonable to conjecture:

\begin{conjecture}
Let $s_i^S(q)$ be the number of sets of $i+1$ distinct integer squares greater than one whose 
sum is below $q$. Similarly let $s_i^C(q)$ be the number of sets of $i+1$ distinct positive integer cubes 
greater than one whose sum is below $q$ then
$$
s_i^S(q) \sim (A_i q)^{\frac{i+1}{2}}
$$
and
$$
s_i^C(q) \sim (B_i q)^{\frac{i+1}{3}}
$$
for constants $1=A_0 > A_1 > \dots > 0$ and $1=B_0 > B_1 > \dots > 0$. 
In both cases the conjecture for $i=0$ holds by the elementary considerations mentioned in the beginning of this section.
\end{conjecture}

\section{Application: Euler characteristics, the M\"obius function and the Riemann Hypothesis}
\label{section: Euler}

In Section~\ref{section: prime complex} we looked at the topology of the Prime complex 
$Prime(q)$ as the quota $q$ was varied. In this section we go back to that complex and consider 
the variance of another topological quantity, the Euler characteristic, of this complex as the quota 
$q$ changes.

Recall for a finite cell complex $X$, the Euler characteristic of $X$, $\chi(X)$, can be defined as 
an alternating sum of the number of cells of various dimensions. Namely, $\chi(X) = \sum_{j=0}^{\infty} (-1)^j c_j \in \mathbb{Z}$
where $c_j$ is the number of $j$-cells in the complex $X$. Recall that the Euler characteristic is a homotopy invariant, i.e., homotopy equivalent spaces have the same Euler characteristic.

Let us first consider, the unrestricted $Prime$ complex where $q=\infty$. This is an infinite dimensional 
simplex, but let us write down an (infinite) expression for its Euler characteristic. When restricted to 
the finite quota complexes $Prime(q)$, we will then get finite sum expressions for their Euler characteristic.

Recall the $j$-dimensional faces of the prime complex are in bijective correspondence with finite sets of $j+1$ distinct primes. Such a finite set of primes can be uniquely associated to its product by the unique 
factorization of integers. Thus the faces of the prime complex are in bijective correspondence with 
the square-free integers $> 1$. (Recall an integer is square-free if no prime is repeated in its factorization).

Recall the M\"obius function $\mu: \mathbb{N} \to \{-1,0,1\}$ is given by
$$
\mu(n)=\begin{cases}
 (-1)^{\text{ number of prime factors of } n }  \text{ if } n \text{ is square-free} \\
 0 \text{ if } n \text{ is not square-free}
\end{cases}
$$

Now note that for a $j$-cell of the prime complex, the corresponding square-free number $n$ 
will have $j+1$ distinct prime factors and so $\mu(n)=(-1)^{j+1}=-(-1)^j$. 
It is now easy to see that formally one has:

$$
\chi(Prime) = \sum_{j=0}^{\infty} (-1)^j c_j = - \sum_{n=2}^{\infty} \mu(n)
$$

To get the corresponding finite sum expression for $\chi(Prime(q))$ we need to include the 
quota restriction which will allow only a finite number of the cells to be in $Prime(q)$.
Well notice a $j$-face $[p_0,p_1,\dots,p_j]$ lies in $Prime(q)$ if and only if 
$\sum_{i=0}^j p_i < q$. Let us define $L_q: \mathbb{N} \to \{0, 1\}$ by 
$$
L_q(n)=\begin{cases}
1 \text{ if } n \text{ is square-free and the sum of its prime divisors is less than } q \\
0 \text{ otherwise}
\end{cases}
$$
Note that $L_q(n)=0$ for sufficiently large $n$ as there are only finitely many primes 
less than $q$ and finitely many square-free numbers with only those prime factors.

Then we have proven:

\begin{pro}[Euler characteristic of the Prime complex]
\label{pro: Eulerprime}
Let $q > 2$ and $Prime(q)$ be the prime complex with quota $q$. 
Then
$$
\chi(Prime(q)) = - \sum_{n=2}^{\infty} \mu(n) L_q(n)
$$
and this sum is finite.
\end{pro}
Note that as $q \to \infty$, $L_q$ converges pointwise to the characteristic function 
of the set of square-free integers. 

%%%%%%Primes: Euler Charactersitic 

\begin{figure}[ht]
\includegraphics[width=0.75\textwidth]{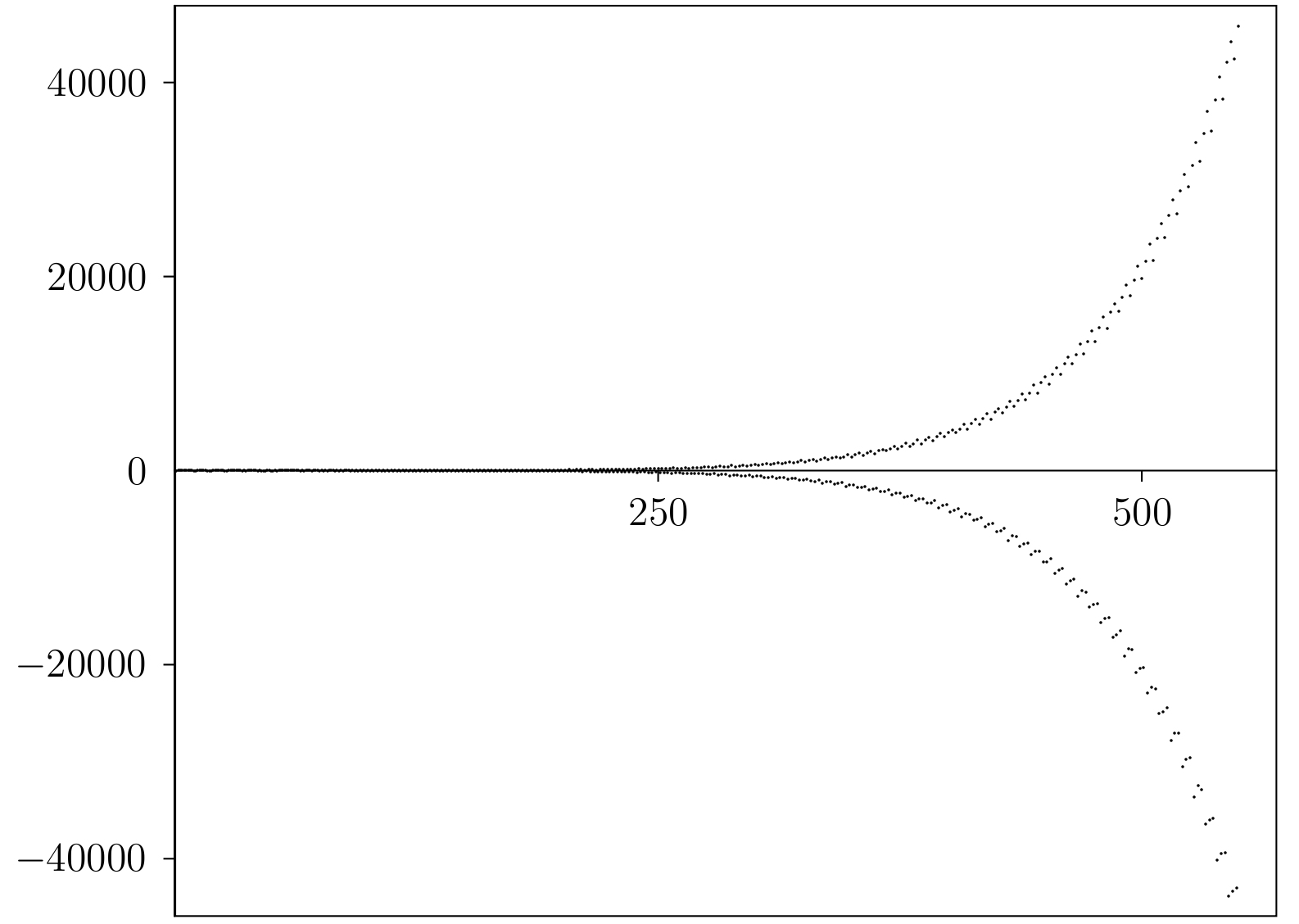}
\caption{$\chi(Prime(q))$, $0<q\le550$}
\label{primesFig:euler}
\end{figure}

In the data the Euler characteristic of $Prime(q)$ behaves very regularly as $q$ increases. 
Figure \ref{primesFig:euler} displays the Euler characteristic for quotas $0<q\le 550$. 
We will see in section~\ref{section: counting} (example 3) that the growth
of $\chi(Prime(q))$ with $q$ is subexponential. 

In analytic number theory, one defines Merten's function
$$
M(N) = \sum_{n=1}^N \mu(n).
$$
By work of Titchmarsh it is known that the Riemann hypothesis that the nontrivial zeros of the 
Riemann zeta function lie on the critical line is equivalent to the statement 
that $M(N) = O(N^{0.5+\epsilon})$ for all $\epsilon > 0$. In other words, for every $\epsilon > 0$, 
there exists a constant $C=C(\epsilon) > 0$ such that 
$$
M(N) \leq C N^{0.5 + \epsilon}
$$
for all $N$ large enough.

With this in mind we introduce a the $LogPrime$ complex which is the full quota complex 
($q = \infty$) with vertices the set of primes where the weight of the vertex $p$ is $ln(p)$. 
Thus as a simplicial complex $LogPrime=Prime$ but $LogPrime(q) \neq Prime(q)$ for most finite 
quotas $q$.

The $j$-dimensional faces of $LogPrime$ are still in bijective correspondence with the 
square-free integers with $j+1$ prime factors. We need only consider the effect of imposing quota 
$q$ with the new weights. Note $[p_0,\dots,p_j]$ is a face of $LogPrime(q)$ if and only if 
$\sum_{i=0}^j ln(p_i) < q$. If $n$ is the square-free number corresponding to the face 
then this is equivalent to $ln(n) < q$ or $n < e^q$.

From this one easily deduces:

\begin{thm}[Euler characteristic of the LogPrime complex]
\label{thm: logprime}
Let $LogPrime(q)$ be the LogPrime complex with quota $q > 2$. Then 
$$
\chi(LogPrime(q)) = - \sum_{2 \leq n < e^q} \mu(n)
$$
If we set $q=ln(N+1)$ then 
$$
\chi(LogPrime(q)) = - \sum_{n=2}^N \mu(n) = 1-M(N)
$$
and so 
$$
|\chi(LogPrime(q))| = O(N^{0.5 + \epsilon})=O(e^{q(0.5 + \epsilon)})
$$
for all $\epsilon > 0$ if and only if the Riemann Hypothesis is true. Equivalently, for any $\epsilon > 0$, 
there is a constant $C=C(\epsilon) > 0$ such that 
$$
ln(|\chi(LogPrime(q))|) \leq (0.5 + \epsilon) q + ln(C(\epsilon))
$$
for all $q$ large enough.
\end{thm}

\begin{figure}[ht]
\includegraphics[width=0.75\textwidth]{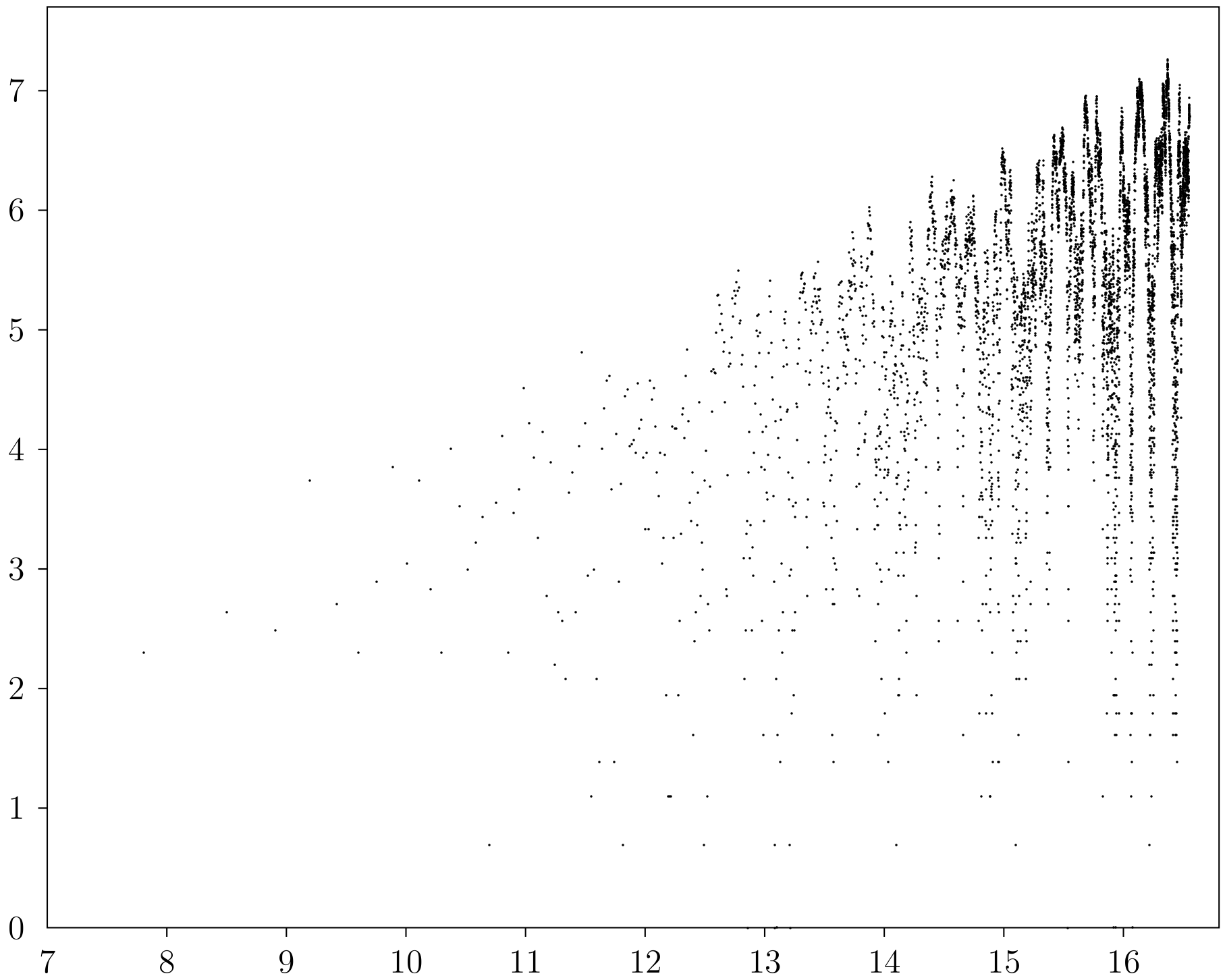}
\caption{$\ln(|\chi(LogPrime(q))|)$, $7\le q\le 16.5554$}
\label{lnprimesFig:euler}
\end{figure}

Figure \ref{lnprimesFig:euler} shows the natural logarithm of the magnitude of the Euler characteristic 
of 6276 of the $LogPrime(q)$ complexes in the range $7\le q\le 16.5554$. 
This data, which was generated using the first $10^6$ odd primes, clearly displays the behavior
expected by the Riemann Hypothesis, namely most data points occur below a line of slope roughly $1/2$. 
Note also that by comparing this data with that for the $Prime$ complex it is clear that $L_q(n)$
has a strongly regularizing effect.

%%%%%%logPrimes: Euler Charactersitic 

From the last theorem, we see that a thorough understanding of how the Euler characteristic 
of the LogPrime complex varies with quota $q$ is equivalent to the Riemann hypothesis.
In fact one direction of the above theorem can be made very explicit. As $\frac{1}{\zeta(s)} = \sum_{n=1}^{\infty} \frac{\mu(n)}{n^s} = 
\sum_{n=1}^{\infty} \frac{M(n)-M(n-1)}{n^s} = 
\sum_{n=1}^{\infty} \frac{\chi(LogPrime(ln(n)))}{n^s} - 
\sum_{n=1}^{\infty} \frac{\chi(LogPrime(ln(n+1)))}{n^s}$
for all $Re(s) > 2$, we see that if the bounds mentioned in Theorem~\ref{thm: logprime} 
hold, the final two L-series would converge for $Re(s) > \frac{1}{2}$ and 
provide an analytic continuation of $\frac{1}{\zeta(s)}$ to $Re(s) > \frac{1}{2}$, hence 
showing that $\zeta(s)$ has no zeros for $Re(s) > \frac{1}{2}$ which is equivalent to the 
Riemann hypothesis.

It would be nice to have a similar characterization of the Riemann Hypothesis using the 
Prime complex itself. As the prime number theorem, twin prime conjecture and Goldbach conjecture 
have topological characterizations using the Prime complex, it iseems likely the Riemann Hypothesis does also but we do not know of any such clean statement currently. However 
section~\ref{section: counting} addresses this further.

\section{Application: The Divisor Complex, perfect, deficient and abundant numbers}
\label{section: divisor complex}

%%%%%%Divisor Graph

\begin{figure}[b]
   \includegraphics[width=\textwidth]{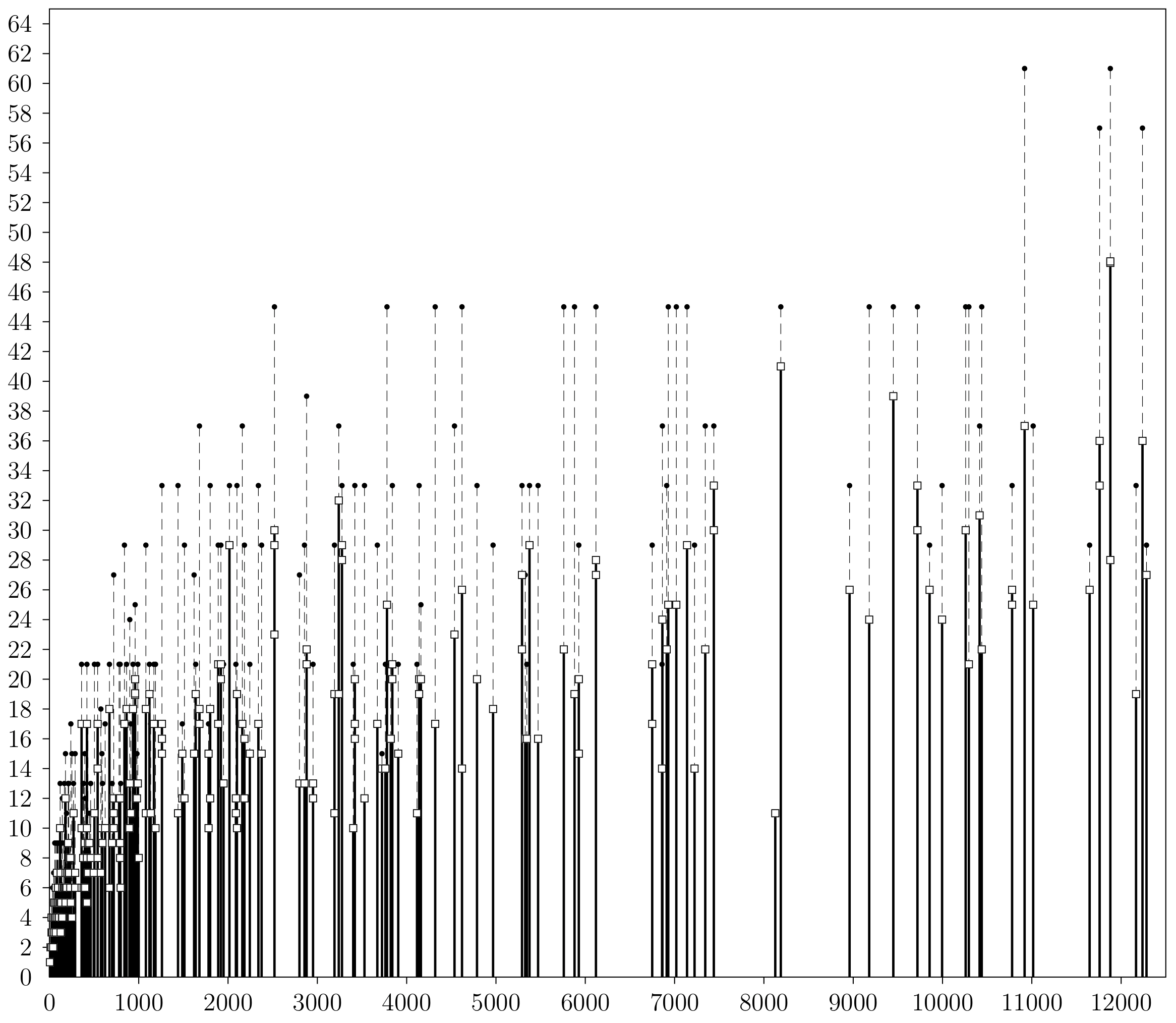}
   \caption{$Div(n)$}
    \label{figure:divisorGraph}
\end{figure}

Recall that a positive integer is perfect/deficient/abundant if the sum of its proper positive divisors 
is equal/less/more than itself.

Consider the divisor complex $Div(n)$ for $n\ge2$. This is the scalar quota complex with 
vertex set the set of proper positive integer divisors of $n$ and quota $n$. So, since the minimal 
vertex has weight 1, $Div(n)$ is homotopic to a bouquet of spheres, consisting of one sphere of 
dimension $i$ for each way $n-1$ can be written as a sum of $i+1$ distinct non-unit divisors of $n$. 
In particular, as stated in the introduction, $n$ is perfect if and only if $Div(n)$ is homotopic to a 
sphere of dimension $\tau(n) -  3$, where $\tau(n)$ is the number of positive integer divisors of $n$. 
Note that $Div(n)$ does not fit into the simple framework developed in section \ref{section: prime complex}
since as the quota $n$ increases the vertex set of $Div(n)$ changes completely. By the Euclid-Euler theorem, 
there is a bijective correspondence between even perfect numbers and Mersenne primes i.e., 
primes of the form $2^p-1$. More specifically a Mersenne prime $2^p-1$ is associated 
to even perfect number $2^{p-1}(2^p-1)$. It is unknown if there are infinitely many Mersenne primes 
(and hence infinitely many even perfect numbers). As of 2010, there were 47 known 
Mersenne primes. It is also unknown if there are any odd perfect 
numbers. 

Figure \ref{figure:divisorGraph} displays $Div(n)$ for $2\le n < 12384$. The data was generated by 
a purpose-built algorithm in C{}\verb!++!. The horizontal axis measures $n$ and the vertical axis measures
topological dimension. Non-contractible complexes are represented as solid vertical lines with 
height equal to the maximal dimension of a sphere in the bouquet; the spheres themselves are plotted
as white squares along the vertical line at height equal to their dimension with spheres of the same 
dimension being listed horizontally. Also, a data point at $(n,\, \tau(n)-3)$ is included for each $n$ 
with $Div(n)$ non-contractible, together with a dashed line connecting it to a highest dimensional sphere in $Div(n)$.
So perfect numbers can be identified as those vertical lines with no dashed component. In this range there are 
four perfect numbers, 6, 28, 496 and 8128, though due to the scale of the graph only 8128 is clearly visible.

One may think of a number $n$ with $Div(n)$ non-contractible as being topologically perfect; 
from this perspective the length of the dashed lines in figure \ref{figure:divisorGraph} measure how close a topologically 
perfect number is to being perfect in the classical sense. The last non-contractible divisor complex 
plotted in figure \ref{figure:divisorGraph} is particularly interesting in this light as it occurs at the odd number 
$n=12285$ and is only a distance of two divisors away from being perfect. A simple check with the data-generating 
algorithm showed that $n=12285$ is the only topologically perfect odd number in the range $2\le n\le 10^6$.
In general of course it is not known if odd perfect numbers exist. 

\section{Combinatorics of Euler characteristics and Lehmer's Conjecture}
\label{section: counting}
In this section, we use the simple combinatorial nature of the Euler characteristic to explicitly 
connect the Euler characteristics of various basic quota systems with well-known combinatorics.
We will work with formal power series but state radii of convergence when reasonable without proof 
(the proofs are elementary).

For the general setup let $1 \leq \nu_1 \leq \nu_2 \leq \dots \leq \nu_i \leq \dots$ be a sequence 
of nondecreasing positive integers with the property that only finitely many of them lie below any 
given $N > 0$.

It is easy to check that we have an identity of formal power series:
$$
\prod_{i=2}^{\infty} (1-x^{\nu_i}) = 1 + \sum_{j=\nu_2}^{\infty} C_j x^j
$$
where 
$$
C_j = \sum_{n=1}^{\infty} (-1)^n D_{n,j} = - \sum_{n=0}^{\infty} (-1)^n D_{n+1,j}
$$
and $D_{n,j}$ is the number of $n$-fold sums of $\nu_i$, $i \neq 1$ with distinct indices, 
which sum to $j$. (We define $C_j=0$ for $ j < \nu_2$ also.)
If $X(q)$ denotes the corresponding scalar quota complex with quota $q > \nu_1$ and $\chi[q]$ its Euler 
chacteristic then by Theorem~\ref{thm: scalarquotathm}, we have
$\chi[q] = 1 - C_{q-1} - \dots - C_{q-\nu_1}$ or equivalently $1-\chi[q] = \sum_{q-\nu_1 \leq i < q} C_i$.
Thus
$$
(1+x+ \dots + x^{\nu_1-1}) \prod_{i=2}^{\infty} (1-x^{\nu_i})
= (1+x+ \dots + x^{\nu_1-1}) + \sum_{j=\nu_2}^{\infty} C_jx^j(1+x+ \dots + x^{\nu_1-1}) 
$$
$$
= (1+x+ \dots + x^{\nu_1-1}) + \sum_{j=\nu_2}^{\infty} (C_j + C_{j-1} + \dots + C_{j-\nu_j+1}) x^j
$$
$$
= (1+x+ \dots + x^{\nu_1-1}) + \sum_{j=\nu_2}^{\infty}  (1-\chi[j+1]) x^j
$$
Multiplying by $1-x$ we get:
$$
\prod_{i=1}^{\infty} (1-x^{\nu_i}) = 1-x^{\nu_1} + (1-x)\sum_{j=\nu_1}^{\infty}(1-\chi[j+1] )x^j
$$
where we used $(1-\chi[j+1])=0$ for $\nu_1 \leq j < \nu_2$. Using the geometric series formula we 
can simplify to get:
$$
\prod_{i=1}^{\infty}(1-x^{\nu_i}) = 1 - (1-x)\sum_{j=\nu_1}^{\infty} \chi[j+1]x^j.
$$
We summarize our results:
\begin{thm}[Generating function for Euler characteristics]
\label{thm: count}
Let $1 \leq \nu_1 \leq \nu_2 \leq \dots \leq \nu_i \leq \dots$ be a sequence of nondecreasing 
positive integers such that only finitely many terms of the sequence are below any given $N > 0$.
Let $\chi[q]$ denote the Euler characteristic of the corresponding scalar quota complex with 
quota $q$. Then $\prod_{i=1}^{\infty} (1-x^{\nu_i})$ defines a well-defined formal power series and 
$$
\prod_{i=1}^{\infty} (1-x^{\nu_i}) = 1- (1-x) \sum_{j=\nu_1}^{\infty} \chi[j+1]x^j
$$
or equivalently
$$
\sum_{j=\nu_1}^{\infty} \chi[j+1]x^j = \frac{1- \prod_{i=1}^{\infty} (1-x^{\nu_i})}{1-x} .
$$
\end{thm}

Theorem~\ref{thm: count} shows that in theory, one can recover the complete quota system from the sequence of Euler characteristics $\chi[q]$. This is because on expanding the product on the left hand side of the first equation, the first nonzero term determines both the lowest weight and its multiplicity by 
comparison with the right hand side of the equation. 
Once this is known, one can divide through by the corresponding factor(s) and determine the 
2nd lowest weight with its multiplicity and proceeding in this manner, recursively recover the full weight system.

\noindent
\textbf{Example 1:} \emph{Counting complex.} Let $\nu_i=i$ for all $i \geq 1$. Then the corresponding 
quota complex will be denoted $Count(q)$. We then  have 
$$
1 - (1-x)\sum_{j=1}^{\infty} \chi(Count(j+1)) x^j = \prod_{n=1}^{\infty} (1-x^n) = \phi(x)
$$
where $\phi(x)$ is called Euler's function (not to be confused with the Euler totient function).
In this case the fomal series and product converge for complex numbers $x$ with $|x| < 1$. 
In fact $\frac{1}{\phi(x)} = \sum_{n=1}^{\infty} p(n) x^n$ where $p(n)$ is the number of partitions 
of $n$. Furthermore writing $x=e^{2 \pi i z}$ makes $\phi(z)$ an analytic function on the upper half 
of the complex plane which turns out to be a modular form. In fact 
$\phi(x) = x^{\frac{-1}{24}} \eta(z)$ where $\eta$ is the Dedekind eta function.\\

\noindent
\textbf{Example 2:} \emph{$Count^{(24)}$ and Lehmer's conjecture.}
Let us take the union of 24 copies of the previous example, 
thus the weights will be the positive integers but there will be 24 vertices for every given weight.
We will denote the corresponding quota complex by $Count^{(24)}(q)$. 
Then
$$
1-(1-x)\sum_{j=1}^{\infty} \chi(Count^{(24)}(j+1))x^j = \prod_{n=1}^{\infty} (1-x^n)^{24} 
= \sum_{n=0}^{\infty} \tau(n+1) x^n
$$
where $\tau(n)$ is Ramanujan's Tau function. From this it is easy to see that Lehmer's conjecture 
that $\tau(n) \neq 0$ for all $n \geq 1$ is equivalent to  
\[\chi(Count^{(24)}(m)) \neq \chi(Count^{(24)}(m+1))\]
for all $m \geq 2$.\\

\noindent
\textbf{Example 3:}  \emph{$Prime$ complex.} Let $\nu_i$ be the $i$th prime. The corresponding quota complex 
is just $Prime(q)$ the prime complex considered in earlier sections. We have 
$$
\sum_{q=2}^{\infty} \chi(Prime(q+1)) x^q = \frac{ 1 - \prod_{p \text{ prime} } (1-x^p)}{1-x} .
$$
The formal power series and product converge for complex numbers $x$ with $|x| < 1$ 
but the sum diverges at $x=1$ and the product diverges at $x=-1$ for example. They define 
analytic functions in the open unit disk and hence on the upper half plane. The fact that 
the radius of convergence is $1$ shows that the growth of $\chi(Prime(q))$ with $q$ is 
subexponential i.e., $\chi(Prime(q))=o(A^q)$ for any $A > 1$.

\section{Random Quota Complexes}
\label{section: Randomquota}

In this section we consider the topology of a random quota complex. 
For the basic facts about probability used in this section see \cite{B} and for the basic 
facts about Fourier transforms see \cite{Ho}.

Fix $X_0=m > 0$ a (nonrandom) value.
Let $X_1, \dots, X_N$ be independent random variables. Furthermore assume $X_i$ are continuous random variables with 
continuous density function $f_i$ with {\bf compact} support in $[m,\infty)$ 
for all $i \geq 1$. Fix a quota $q  > m$. Since we assume the densities are continuous with compact support, $f_i \in L^1(\mathbb{R}) \cap L^2(\mathbb{R})$ for $1 \leq i \leq N$ and we will hence 
be able to use convolutions, Fourier transforms and the inverse Fourier transform without technical 
difficulties. We will do so without further mention. Note one can consider the density for 
$X_0=m$ as the delta measure centered at $m$ but while this has a well-defined Fourier transform, 
the inverse transform formula is singular and we will worry about this when necessary.

\begin{defn}[Random scalar quota complex]
Let $\mathbb{X}=\{ X_0=m,X_1,\dots,X_N \}$ be chosen as in the previous paragraph and 
$q > m > 0$. 
$\mathbb{X}[q]$ is the quota complex on vertices $\{ 0, 1, 2, \dots, N \}$ with 
weights $w(i)=X_i$ and quota $q$.

$\mathbb{X}[q]$ is called a random scalar quota complex. On any run of the experiment, the weights will 
take on specific values and $\mathbb{X}[q]$ will determine a specific scalar quota complex. 
Thus $\mathbb{X}[q]$ can be viewed as a random variable on the sample space of the underlying experiment which takes values in the set of finite abstract simplicial complexes.
\end{defn}

By Theorem~\ref{thm: scalarquotathm}, we see immediately that $\dim(\bar{H}_j(\mathbb{X}[q],\mathbb{Q}))$ is an integer valued random variable whose value is the number of $(j+1)$-fold 
sums of the variables $X_1,\dots, X_N$ (repeats not allowed) that lie in the interval 
$[q-m,q)$. The number of $j$-dimensional faces of $\mathbb{X}[q]$ is another integer random variable 
whose value is the number of $(j+1)$-fold sums of the variables $X_0, \dots, X_N$ (repeats not allowed) 
that are below quota $q$. Note that as each random variable takes values greater than or equal to $m > 0$, there 
will be no $j$-faces when $q \le (j+1)m$ and so the dimension of $\mathbb{X}[q]$ is bounded 
by $\frac{q}{m}-1$.

Let $\mathfrak{J} \subseteq \{ 1, \dots, N \}$ have $|\mathfrak{J}|=j$ and let 
$X_\mathfrak{J} = \sum_{i \in \mathfrak{J}} X_i$ be the corresponding $j$-fold sum. 

Let us form a Bernoulli indicator random variable:
$$
B_{\mathfrak{J}} = 
\begin{cases}
1 \text{ if } X_{\mathfrak{J}} \in [q-m,q) \\
0 \text{ otherwise}
\end{cases}
$$

Then calculating expected values we see:
$$
E[B_{\mathfrak{J}}] = Pr(X_{\mathfrak{J}} \in [q-m,q) ) 
$$
and so to figure out this expectation, we need to determine the distribution for the corresponding 
sum $X_{\mathfrak{J}}$. This is standard but we include a quick exposition here:

\begin{defn} Let $f, g \in L^1(\mathbb{R})$ then we define the convolution 
$$
( f \star g )(\alpha) = \int_{-\infty}^{\infty} f(\alpha - x) g(x) dx
$$
for all $\alpha \in \mathbb{R}$. It is well known that $f \star g \in L^1(\mathbb{R})$ and indeed 
$ || f \star g ||_1 \leq || f ||_1 || g ||_1 $. $(L^1(\mathbb{R}), \star)$ forms a commutative 
Banach algebra.
\end{defn}

\begin{defn} Let $f$ be a continuous density function for a real random variable $X$. 
The cumulative density function $F$ is defined as 
$$
F(\alpha) = Pr( X \leq \alpha) = \int_{-\infty}^{\alpha} f(x)dx
$$
and so we have $F'(\alpha)=f(\alpha)$ for all $\alpha \in \mathbb{R}$.
\end{defn}

 The relevant proposition is:
 
 \begin{pro}
 \label{pro: conv}
 Let $X$ and $Y$ be {\bf independent} continuous real valued random variables with continuous density 
 functions $f_X, f_Y$ respectively and cumulative density functions $F_X, F_Y$ respectively.
 Then $F_{X+Y} = F_X \star f_Y$ and $f_{X+Y} = f_X \star f_Y$.
 \end{pro}
 \begin{proof}
 Computing we find:

 \begin{align*}
 \begin{split}
 F_{X+Y}(\alpha) &= \int \int_{ x+y \leq \alpha} f_X(x)f_Y(y) dx dy \\
 &= \int_{-\infty}^{\infty} \left(\int_{-\infty}^{\alpha - y} f_X(x)dx \right) f_Y(y) dy  \\
 &= \int_{-\infty}^{\infty} F_X(\alpha-y) f_Y(y) \\
 &= (F_X \star f_Y) (\alpha)
 \end{split}
 \end{align*}
 The corresponding identity for $f_{X+Y}(\alpha)=F_{X+Y}'(\alpha)$ is obtained by 
 differentiating the 2nd row of equations above with respect to $\alpha$.
 \end{proof}

\begin{cor}
Let $X_1, \dots X_N$ be independent real random variables with continuous densities 
$f_1, \dots f_N$. Then $f_{X_1+\dots + X_N} = f_1 \star f_2 \star \dots \star f_N$
\end{cor}
\begin{proof}
Follows by induction on Proposition~\ref{pro: conv} by setting $X=X_1$ and 
$Y=X_2 + \dots + X_N$.
\end{proof}

Returning to the calculation of $E[B_{\mathfrak{J}}]=Pr(X_{\mathfrak{J}} \in [q-m,q))$ 
let us define $f_{\mathfrak{J}}$ to be the convolution of $f_i$ for $i \in \mathfrak{J}$. 
Thus $f_{\mathfrak{J}}$ is a $|\mathfrak{J}|$-fold convolution and 
$E[B_{\mathfrak{J}}] = \int_{q-m}^q f_{\mathfrak{J}} dx$.

Since 
$$
\dim(\bar{H}_{j-1}(\mathbb{X}[q])) = \sum_{\mathfrak{J} , |\mathfrak{J}|=j} B_{\mathfrak{J}}
$$
we find 
$
E[\dim(\bar{H}_{j-1}(\mathbb{X}[q], \mathbb{Q}))] = \sum_{\mathfrak{J}, | \mathfrak{J}|=j } \int_{q-m}^q f_{\mathfrak{J}} dx 
$. 

Let 
$$
\mathbb{I}_m (x) = \begin{cases} 
1 \text{ if }  x \in [0,m] \\
0 \text{ otherwise}
\end{cases}
$$
Then it is easy to check that $(f \star \mathbb{I}_m)(q) = \int_{q-m}^q f(x) dx$.
Thus we conclude 
$$
E[\dim(\bar{H}_{j-1}(\mathbb{X}[q], \mathbb{Q}))] = \sum_{\mathfrak{J}, | \mathfrak{J}|=j } (f_{\mathfrak{J}} 
\star \mathbb{I}_m) (q).
$$
Since the convolution of continuous functions with compact support is continuous with compact 
support, we see that $E[\dim(\bar{H}_{j-1}(\mathbb{X}[q], \mathbb{Q}))]$ is a continuous 
function of $q$ with compact support.  (Though $\mathbb{I}_m$ is not continuous, one can 
go back to the integral expression a couple of lines back to check the continuity of the 
expected dimension of homology.)

Since the Euler characteristic can be expressed as the alternating sum of these homology 
groups, we have 
$$
E[\chi(\mathbb{X}[q])] = 1+\sum_{\mathfrak{J} \subseteq \{1, \dots, N \} } (-1)^{|\mathfrak{J}|-1} (f_{\mathfrak{J}}\star \mathbb{I}_m)(q)
$$
is a continuous function of $q$ with compact support.

Note we do not include the empty set $\mathfrak{J}$ in the sum and the extra $1$ in front on the 
right hand side is due to 
$\dim(H_0)=\dim(\bar{H}_0) + 1$, or equivalently to account for the vertex of minimal weight $m$.

Rewriting a bit we get:
$$
1-E[\chi(\mathbb{X}[q])] = \sum_{\mathfrak{J} \subseteq \{1, \dots, N \} } (-1)^{|\mathfrak{J}|} (f_{\mathfrak{J}}\star \mathbb{I}_m)(q).
$$
Note that though the left hand side is defined only for $q > m$, the right hand side is 
defined for all $q$ and so can be viewed as a continuation of the left hand side to the whole real line.
To shed light we will take the Fourier transform of the function
$$
G(q) = \sum_{\mathfrak{J} \subseteq \{1, \dots, N \} } (-1)^{|\mathfrak{J}|} (f_{\mathfrak{J}}\star \mathbb{I}_m)(q)
$$
appearing on the right hand side. (Note that Laplace transforms would work just as well for our purposes.)

As the fourier transform of a convolution is the product of 
the fourier transforms, we have
$$
\hat{G}(\alpha) = \sum_{\mathfrak{J} \subseteq \{1, \dots, N \} } (-1)^{|\mathfrak{J}|} 
\hat{f_{\mathfrak{J}}}(\alpha) \hat{\mathbb{I}}_m(\alpha)
$$
However now we have $\hat{f}_{\mathfrak{J}}(\alpha) = \prod_{i \in \mathfrak{J}} \hat{f}_i$. 
Using this it is easy to see that 
$$
\hat{G}(\alpha) = \left[\prod_{j=1}^N (1 - \hat{f}_j) - 1 \right] \hat{\mathbb{I}}_m(\alpha)
$$
where 
$\hat{\mathbb{I}}_m(\alpha) = \frac{1-e^{-2 \pi i \alpha m}}{2\pi i \alpha}$ 
is found by direct computation. To make the formula more uniform, we note that 
$X_0=m$ can be thought of as having density function given by the delta measure 
centered at $m$, i.e., $f_0=\delta_m$. The Fourier transform of this measure, $\hat{f}_0(\alpha)$,  
is given by $\int_{-\infty}^{\infty} \delta_m(x)e^{-2\pi i \alpha x} dx = e^{-2 \pi i \alpha m}$.
Thus we can write $\hat{I}_m = \frac{1 - \hat{f}_0}{2 \pi i \alpha}$ and so
$$
\hat{G}(\alpha) = \frac{1}{2\pi i \alpha}\left[\prod_{j=0}^N (1-\hat{f}_j(\alpha)) - (1-\hat{f}_0(\alpha))\right].
$$
(The reader is warned that though $\hat{G}(\alpha)$ has a continuous Fourier inverse with compact 
support, $(1-\hat{f}_0)$ by itself has badly 
behaved inverse Fourier transform.) By the Fourier inversion formula we get:
$$
1 - E[\chi(\mathbb{X}[q])] = \frac{1}{4 \pi^2 i} \int_{-\infty}^{\infty} e^{2 \pi i \alpha x} \left(\prod_{j=0}^N (1-\hat{f}_j(\alpha)) - (1-\hat{f}_0(\alpha))\right) \frac{d\alpha}{\alpha}.
$$

We will see in an example that this is an example of an Euler product decomposition 
in analytic number theory but first let us record the probability results we have obtained.

\begin{thm}[Expected topology of random scalar quota complexes]
\label{thm: ExpTop}
Let $X_0=m > 0$. Let $X_1, \dots, X_N$ be independent, continuous random variables 
with density functions $f_1, \dots, f_N$ which are continuous with 
compact support in $[m, \infty)$ and let $\mathbb{X}[q]$ be the random scalar quota 
complex determined by this collection and quota $q > m > 0$.

Then for $j \geq 1$, 
$$
E[\dim(\bar{H}_{j-1}(\mathbb{X}[q], \mathbb{Q}))] = \sum_{\mathfrak{J}, |\mathfrak{J}|=j} (f_{\mathfrak{J}} \star \mathbb{I}_m)(q)
$$
is a continuous function of $q$ with compact support.

Furthermore we have 
$$
1 - E[\chi(\mathbb{X}[q])] = \frac{1}{4 \pi^2 i} \int_{-\infty}^{\infty} e^{2 \pi i \alpha x} \left(\prod_{j=0}^N (1-\hat{f}_j(\alpha)) - (1-\hat{f}_0(\alpha))\right) \frac{d\alpha}{\alpha}.
$$
is a continuous function of $q$ with compact support.

\end{thm}

We now present a degenerate example to illustrate that in some sense in Theorem~\ref{thm: ExpTop}, 
the last equality has right hand side a form of Euler product.

First recall (see \cite{A}), the definition of the Riemann zeta function and its reciprocal:
$$
\zeta(s) = \sum_{n=1}^{\infty} \frac{1}{n^s}
$$
converges for all complex numbers $s$ with $Re(s) > 1$. This is the L-series associated to the 
constant function with value 1 (as the numerators are all 1's). It is well known and easy to check that 
due to the unique factorization of positive integers into primes, we have an Euler product for 
the zeta function:
$$
\zeta(s) = \sum_{n=1}^{\infty} \frac{1}{n^s} = \prod_{p \in Prime} \frac{1}{1-\frac{1}{p^s}}
$$
which also holds for $Re(s) > 1$.
Similarly, the reciprocal of the Riemann zeta function is the L-series associated to the M\"obius function 
and so
$$
\psi(s)=\frac{1}{\zeta(s)} = \sum_{n=1}^{\infty} \frac{\mu(n)}{n^s} = \prod_{p \in Prime} \left(1-\frac{1}{p^s}\right)
$$
for $Re(s) > 1$. It is common to look at restricted Euler products to various sets of primes and 
we will consider the simple case where $P$ is a finite set of primes. Let us write 
$$
\psi_P(s) = \prod_{p \in P} \left(1-\frac{1}{p^s}\right).
$$
Note $\psi_P$ is just a finite product of entire functions so is entire with zeros only lying on the 
imaginary axis of the form $\frac{2 \pi i k}{ln(p)}$ where $p \in P$ and $k$ an integer.

Consider the case where $P$ is the set of primes $\leq N$ for some $N > 3$ and set quota $q=ln(N+1)$. Consider 
``random" variables $X_0=ln(2), X_1=ln(3), X_2=ln(5), \dots$ which are just constant at the log-prime values.
The corresponding density functions are delta measures centered at the corresponding log-primes 
and the corresponding ``random" complex is just the LogPrime complex. 
We will apply Theorem~\ref{thm: ExpTop} to this scenario. If the reader is worried about the 
the fact that the density functions $f_i$ are not continuous with compact support, just replace 
the delta measures with very small continuous bumps of mass 1 and compact support 
localized around these log-primes. Since there are no serious benefits of the more rigorous approach we will just use the delta measures for this example. 

By Theorem~\ref{thm: logprime}, we have 
$1-\chi(LogPrime(q)) = \sum_{1 \leq n < e^q} \mu(n)$.
Note the Fourier transform of the delta measure $\delta_{ln(p)}$ centered at $ln(p)$ is 
$e^{-2 \pi i ln(p) \alpha}=p^{-2 \pi i \alpha}$ and so the formula from Theorem~\ref{thm: ExpTop} becomes
\begin{align*}
\begin{split}
\sum_{1 \leq n < e^q} \mu(n) = \frac{1}{4 \pi^2 i} \int_{-\infty}^{\infty} e^{2 \pi i \alpha x} 
\left(\prod_{p \in P}(1-p^{-2 \pi i \alpha}) - (1-2^{-2 \pi i \alpha})\right) \frac{d \alpha}{\alpha} .
\end{split}
\end{align*}
Doing a change of variable $s=2 \pi i \alpha$ one finds:
$$
\sum_{1 \leq n < e^q} \mu(n) = \frac{1}{4 \pi^2 i} \int_{-i\infty}^{+i\infty} e^{sx}\left(\psi_P(s) - \left(1-\frac{1}{2^s}\right)\right) 
\frac{ds}{s}.
$$

\section{Concluding Remarks}

We have seen the standard Morse theoretic problem of studying the change of topology 
of a space as a parameter is varied can be applied to the study of quota complexes as the 
quota is varied. In this case the quota determines the complex through a set of linear inequalities.
A similar though more difficult nonlinear version of this problem is encountered in the study 
of random data sets through the method of persistant homology. 

In this context a finite set of data points in $\mathbb{R}^n$ is studied by considering the 
space $X[r]$ consisting of the union of open balls of radius $r$ around these points. The topology 
of $X[r]$ changes from being discrete when $r$ is very small to a single contractible blob when $r$ is very large. 
Usually $r$ is varied in a range corresponding to reasonable error estimates for the experiment and 
the persistant features of the topology of $X[r]$ as $r$ ranges over this error interval are looked for 
and attributed as features of the data set.

In this context one often forms the associated Vietoris-Rips complex instead of the actual 
space $X[r]$ for computational simplicity. This is an abstract simplicial complex 
where there is one vertex for each of the original data points and where 
$[v_0,\dots,v_n]$ is a face if and only if the radius $r$ balls around the $v_i$ have pairwise nonempty intersection. This is equivalent to $D(v_0,\dots,v_n)=sup_{i,j} |v_i - v_j| < 2r$. Thus instead 
of the linear inequalities present in scalar quota complexes, these inequalities are nonlinear 
ones saying the diameters of various sets should be below quota $q=2r$.
In these cases, one has to do more work to study the change of topology 
as for example Theorem~\ref{thm: scalarquotathm} does not apply. 

Please see \cite{Ka} for more details. We will not pursue this here.

\appendix
\section{Vector-weighted quota complexes.}
\label{section: vector}

In order to state topological structure results about vector weighted quota complexes, 
we need to recall the definition of the Lusternik-Schnirelmann category of a space.

\begin{defn}
Let $X$ be a topological space. Let $n$ be a nonnegative integer. 
We say that $Cat(X) \leq n$ if $X$ can be covered by $n+1$ open sets 
$U_0, \dots, U_n$ such that each $U_i$ is contractible in $X$, i.e., the 
inclusion map $j_i: U_i \to X$ is homotopically trivial.

The category of a space is hence either a nonnegative integer or infinite.
Note $Cat(X)=0$ if and only if $X$ is contractible. $Cat(X)$ is a homotopy invariant 
of $X$ i.e., homotopy equivalent spaces have the same category.
\end{defn}

The notion of category was introduced to give a lower bound on the number of critical points 
of a Morse function on a smooth manifold. It has also proved fruitful in proving bounds on 
degrees of nilpotence of various algebraic structures associated to a space $X$.
(See \cite{PY} for example.)

Note any closed $n$-manifold can be covered by a finite number of open sets homeomorphic 
to $\mathbb{R}^n$ and hence contractible. Thus all closed manifolds have finite category.
A sphere $S^n$ has $Cat(S^n)=1$ as a sphere can always be covered by two contractible open 
sets (an upper and lower hemisphere) but is itself not contractible. If $\mathbb{C}P^n$ is 
$n$-dimensional complex projective space then $\mathbb{C}P^n$ can be covered 
by $n+1$ charts which are contractible pieces so $Cat(\mathbb{C}P^n) \leq n$. 
However cup product arguments can be used to show that $Cat(\mathbb{C}P^n) > n-1$ 
and so $Cat(\mathbb{C}P^n)=n$ for all integers $n \geq 0$.

A bouquet of positive dimensional 
spheres $X$ has $Cat(X)=0$ if it is an empty bouquet (just a point) or 
$Cat(X)=1$ otherwise. To see this let $U_0$ be an open thickening of the attaching point 
of the bouquet and $U_1$ be $X$ minus the attaching point. $U_0$ is contractible 
and though $U_1$ is not in general contractible, it is always contractible in $X$ i.e., the 
inclusion map $U_1 \to X$ is homotopically trivial. This is because each component of 
$U_1$ is contractible and $X$ is path connected as the bouquet does not involve $0$-spheres.

A bouquet with $0$-spheres involved can have category larger than $1$ as the bouquet 
will not be path connected. Each $0$ sphere gives an extra point component which increases 
the category by $1$.  To address the issue of $0$-spheres in quota complexes we 
define the concept of a "shell vertex" in a quota complex.

\begin{defn}[Shell vertices]
Let $X$ be a finite scalar quota complex with vertex of minimal weight $v_{min}$.
A vertex $s$ of $X$ is called a shell vertex if 
$q-w(v_{min}) \leq w(s) < q$. 

Note by Theorem~\ref{thm: scalarquotathm}, $X$ is homotopy equivalent to a bouquet 
of spheres where there is one $0$-sphere for each shell vertex not equal to $v_{min}$. Thus 
$Cat(X)=(\text{Number of shell vertices not equal to } v_{min}) + \epsilon$ where $\epsilon=1$ 
if there are any positive dimensional spheres in the bouquet and $\epsilon=0$ if not.

For a vector valued quota complex, for the $i$th coordinate of the weights and quota we 
can find a vertex $v_{min,i}$ of minimal weight for that coordinate. Then a vertex $s$ is called a 
shell vertex if $q_i-w_i(v_{min,i}) < w_i(s) < q_i$ for {\bf some} coordinate $i$.

\end{defn}

\begin{thm}[Category of vector weighted quota systems]
Let $X$ be a finite vector weighted quota system with weight dimension $N$, weight function $\hat{w}$ and quota $\hat{q}$.
 So $\hat{w}: V \to \mathbb{R}_+^N$ and quota $\hat{q} \in \mathbb{R}_+^N$.
 
 $X$ is then the union of $N$ scalar quota complexes.
 
 Furthermore 
 if $X$ has no shell vertices, we have $Cat(X) \leq 2N-1$. Thus $\frac{Cat(X)+1}{2}$ provides a 
 homotopy invariant lower bound for the weight dimension of a quota complex with no shell 
 vertices.
 
\end{thm}
\begin{proof}
Let $X_i$ be the complex determined by the $i$th component of the weight-quota system 
for $1 \leq i \leq N$. Then by definition $X = \cup_{i=1}^N X_i$. By Theorem~\ref{thm: scalarquotathm}, 
each $X_i$ is homotopy equivalent to a bouquet of spheres. If $X$ has no shell vertices then no 
$0$-spheres occur in these bouquets and hence each $X_i$ can be covered by two open 
sets which are contractible in $X_i$ and hence in $X$. Thickening up these open sets so that they are 
open in $X$, we see that $X$ can then be covered by $2N$ open sets, each contractible 
in $X$. Thus $Cat(X) \leq 2N-1$ by definition. The rest of the theorem then follows readily.

\end{proof}

\section{Every finite simplicial complex is a quota complex.}
\label{section: Alves}

In this appendix we provide a proof that every finite simplicial complex is a quota complex.
This fact already appeared in unpublished work of Manuel Alves who was working on independent undergraduate 
research on voting theory under the guidance of the first author. Given a monotone 
voting system on a set $V$ of voters, the voting complex associated to the system 
is an abstract simplicial complex based on $V$ whose faces consist of the losing coalitions 
of the system.

The proof of the following theorem then parallels the fact that every monotone voting system 
is a (vector valued) quota system. Let $|F|$ denote the number of vertices in a face $F$.

\begin{thm}
If $X$ is a finite simplicial complex on vertex set $V$ then $X$ is a quota complex i.e., 
there exists a weight function $\hat{w}: V \to \mathbb{R}^n_+$ for some $n$ and quota 
$\hat{q} \in \mathbb{R}^n_+$ such that $X$ is the simplicial complex associated to 
quota system $[\hat{w}: \hat{q}]$.

Furthermore the weights and the quota 
can be chosen to have positive integer coordinates and these coordinates can be assumed 
to all be distinct.

\end{thm}
\begin{proof}
Let $\mathfrak{L}$ denote the set of maximal faces (facets) of $X$. 
Order these faces $\mathfrak{L}=\{F_1, \dots, F_s \}$. 
Define a scalar valued quota system 
$w_i: V \to \mathbb{R}_+$ by setting 
$$
w_i(v)=\begin{cases}
1 \text{ if } v \in F_i \\
|F_i|+1 \text{ if } v \notin F_i
\end{cases}
$$
Set quota $q_i=|F_i|+1$.
It is simple to check that the subcomplex of the simplex on $V$ determined by the quota system 
$[w_i : q_i]$ is exactly the simplex $F_i$ and its faces. Letting $\hat{w}: V \to \mathbb{R}_+^s$ 
be the function with coordinates the $w_i$ and $\hat{q} \in \mathbb{R}_+^s$ be the 
vector with coordinates the $q_i$, then $[\hat{w} : \hat{q}]$ is subcomplex of the simplex 
on $V$ given by $\bigcup_{i=1}^s F_i$. Hence $X=\bigcup_{i=1}^s F_i$ is the quota complex 
$[\hat{w}:\hat{q}]$.

Define as usual the weight of a face to be the sum of the weights of the vertices in the face.
Since $X$ has a finite number of faces, only a finite number of vectors occur as 
weights of the faces of $X$. It is clear then that we may 
change the coordinates of the quota by subtracting sufficiently small positive numbers from the coordinates 
so that the quota determines the same complex as the original quota system and has distinct positive rational entries and that these quota coordinates are not 
acheived as the weight of any face in the simplex on $V$. It is then clear that we can perturb 
the weights of the vertices by small amounts so that the resulting quota complex still 
determines the same complex and so that the weights of the vertices consist of 
vectors of distinct positive rational numbers and such that no weight vectors share 
any coordinates which each other or with the quota vector.

Finally scaling all weight vectors and the quota vector by the least common multiple of 
the denominators of all rational numbers involved, we get a quota system determining the same 
complex as the original, where weights and 
quota have positive integer coordinates and all these coordinates are distinct. 

\end{proof}

\section{Expectation Justification}
\label{section: ExpJust}

We adopt the setup and notation from the beginning of section~\ref{section:  prime complex} in this section of the appendix.

Furthermore set  
$S_i^{high}(x) = \widehat{s}_i(x)/2^{f(x/(i+1))}$ and $S_i^{low}(x) = \widehat{s}_i(x/(i+1))/2^{f(x)}$. We have that 
\[ S_i^{low}(q) \approx \frac{ {s_0(\lceil q/(i+1) \rceil) \choose i+1}}{\sum_j  {s_0(q) \choose j}} \le S_i(q) \le \frac{ {s_0(q) \choose i+1}}{\sum_j  {s_0(\lceil q/(i+1) \rceil) \choose j}} \approx S_i^{high}(q) .\]
So if we take $\widehat{S}_i(x) = {f(x) \choose i+1} / 2^{f(x)}$ then for $x$ with $f(x/(i+1)) \ge i$, $S_i^{low}(x) < \widehat{S}_i(x) < S_i^{high}(x)$.  
The global behavior of the family of functions $\widehat{S}_i(x)$ is easy to ascertain; we collect some of the features in the 
following proposition.

\begin{pro}
\label{prime complex: proposition}
Let $k'$ be the smallest positive integer such that $f(\kappa) \le k'$, and for $j\ge k'$ let $x_j\in[\kappa,\, \infty)$ 
such that $f(x_j) = j$. Then for all $i \ge k'$, $\widehat{S}_i(x_{i}) = 0$, $\lim_{x\to \infty} \widehat{S}_i(x) = 0$ and 
$\widehat{S}_i(x)$ has exactly one critical value in $[x_{i},\,\infty)$ which is a maximum. Moreover, if $m_i$ is  the 
critical point of $\widehat{S}_i(x)$ in $[x_{i},\,\infty)$, then $x_{2i+1} < m_i < x_{2i+2}$ and $\lim_{i\to \infty}\widehat{S}_i(m_i) = 0$.
\end{pro}

\begin{proof} 
The first assertions are an exercise in undergraduate calculus: we may compute
\[ \widehat{S}_i'(x) = f'(x)\widehat{S}_i(x)\left[ \sum_{j=0}^{i} \frac{1}{f(x)-j} -\ln 2 \right], \]
from which it follows that $\widehat{S}_i(x)$ has exactly one critical value in $[x_{i},\,\infty)$ which is a maximum. 

In order to establish the bounds on $m_i$ we use the following sequence (see \cite{BBG} pg. 10)
\[ l_i = \sum_{j=2^i}^{2^{i+1}-1} \frac{1}{j} \to \ln2 \hspace{0.1in}\text{as}\hspace{0.1in}i\to\infty.\]
Namely, consider the sequences 
\[a_i = \sum_{j=0}^{i}\frac{1}{f(x_{2i+1})-j} = \sum_{j=i+1}^{2i+1}\frac{1}{j} \hspace{0.1in}\text{and}\hspace{0.1in} b_i = \sum_{j=0}^{i}\frac{1}{f(x_{2i+2})-j} = \sum_{j=i+2}^{2i+2}\frac{1}{j}. \]
Then, $a_i$ is a bounded decreasing sequence and $b_i$ is a bounded increasing sequence, so they converge. 
Furthermore $a_{2^i-1} = l_i$ and $b_{2^i-2} - l_i \to 0$ as $i\to \infty$, so $a_i, b_i \to \ln2$ as $i\to\infty$. Hence, 
since $a_i$ is decreasing we must have that $x_{2i+1} < m_i$, and since $b_i$ is increasing we must have that 
$m_i<x_{2i+2}$. So $\widehat{S}_i(m_i) < {2i+2 \choose i+1}/2^{2i+1}$ and since, by Stirling's formula, 
${2n \choose n} \sim 2^{2n}/\sqrt{\pi n}$ it follows that $\lim_{i\to\infty}\widehat{S}_i(m_i) = 0$.
\end{proof}

Note that while it is reasonable to make the the approximation
\[ H_i(q) =  \frac{s_i(q) - s_i(q-v_1)}{\sum_j s_j(q) - \sum_j s_j(q-v_1)}  \approx \frac{{f(q) \choose i+1} - {f(q-v_1) \choose i+1}}{2^{f(q)} - 2^{f(q-v_1)}},\]
the function on the right seems to be less amenable to a simple general analysis.

\bigskip

\noindent
Dept. of Mathematics \\
University of Rochester, \\
Rochester, NY 14627 U.S.A. \\
E-mail address: jonpak@math.rochester.edu \\

\bigskip

\noindent
Dept. of Mathematics \\
University of Rochester, \\
Rochester, NY 14627 U.S.A. \\
E-mail address: winfree@math.rochester.edu \\


\begin{thebibliography}{EMG}

\bibitem[A]{A}{T. M. Apostol}
\textit{Introduction to Analytic Number Theory, Undergraduate Texts in Mathematics,}
New York-Heidelberg: Springer-Verlag, 1976.

\bibitem[B]{B}{P. Billingsley}
\textit{Probability and Measure,}
Wiley, New York, 1986.

\bibitem[BBG]{BBG}{J. Borwein, D. Bailey, R. Girgensohn}
\textit{Experimentation in Mathematics: Computational Paths to Discovery,}
A K Peters, Wellesley, MA, 2004.

\bibitem[EH]{EH}{H. Edelsbrunner, J. Harer,}
\textit{Persistant Homology -  a Survey.}
In \textit{Surveys on discrete and computational geometry}, volume 453 of \textit{Contemp. Math}, pages 257-282.
Amer. Math. Soc., Providence RI, 2008.   

\bibitem[H]{H}{A. Hatcher}
\textit{Algebraic Topology}
Cambridge University Press, Cambridge 2002.

\bibitem[Ho]{Ho}{K. B. Howell}
\textit{Principles of Fourier Analysis}
Chapman \& Hall/CRC, Boca Raton, FL, 2001.

\bibitem[Ka]{Ka}{M. Kahle}
\textit{Random Geometric Complexes}
to appear in Discrete \& Computational Geometry.

\bibitem[LR]{LR}{G. Latouche, V. Ramaswami}
\textit{Introduction to Matrix Analytic Methods in Stochastic Modelling,} 1st ed. 
Chapter 1: Quasi-Birth-and-Death Processes; ASA SIAM, 1999.

\bibitem[Mu]{Mu}{J. R. Munkres}
\textit{Elements of algebraic topology}
Addison-Wesley, Menlo Park, CA, 1984.

\bibitem[PY]{PY}{J. Pakianathan, E. Yalcin,}
\textit{On nilpotent ideals in the cohomology ring of a finite group,} 
Topology {\bf 42}, (2003), 1155-1183.

\bibitem[Ti]{Ti}{E .C. Titchmarsh}
\textit{The theory of the Riemann zeta-function,} 2nd ed. revised by D.R. Heath-Brown,
Oxford University Press, New York,1986.

\bibitem[Te]{Te}{G. Tenenbaum}
\textit{Introduction to analytic and probabilistic number theory,}
Cambridge University Press, Cambridge 1995.

\end{thebibliography}
\end{document}